\documentclass[bj]{imsart}

\arxiv{0000.0000}
\usepackage{amsthm}
\usepackage{amsmath}
\usepackage{amssymb}
\usepackage{natbib}
\usepackage[colorlinks,citecolor=Blue,urlcolor=blue,filecolor=blue,linkcolor=red,backref=page]{hyperref}
\RequirePackage{accents}
\RequirePackage{amsfonts,bm}
\usepackage[usenames,dvipsnames,svgnames,table]{xcolor}
\usepackage{booktabs}
\usepackage{bbm, dsfont}
\usepackage{mathrsfs} 
\makeatletter
\renewcommand*{\@fnsymbol}[1]{\ensuremath{\ifcase#1\or \dagger\or \dagger\or \ddagger\or
    \mathsection\or \mathparagraph\or \|\or **\or \dagger\dagger
    \or \ddagger\ddagger \else\@ctrerr\fi}}
\makeatother

\newtheorem{theorem}{Theorem}[section]
\newtheorem{corollary}[theorem]{Corollary}
\newtheorem{lemma}[theorem]{Lemma}
\newtheorem{assumption}{Assumption}
\newtheorem{definition}[theorem]{Definition}

\newtheorem{remark}{Remark}

\DeclareMathOperator{\Tr}{Tr}

\newcommand{\bbeta}{\boldsymbol \beta}

\newcommand{\bGamma}{\boldsymbol \Gamma}
\newcommand{\bSigma}{\boldsymbol \Sigma}
\newcommand{\bOmega}{\boldsymbol \Omega}
\newcommand{\bX}{\boldsymbol X}

\newcommand{\eig}{\mathrm{eig}}

\newcommand{\bPhi}{\boldsymbol \Phi}
\newcommand{\bPsi}{\boldsymbol \Psi}

\newcommand{\bI}{\boldsymbol I}
\newcommand{\vc}{\text{Vec}}
\newcommand{\bDel}{\boldsymbol \Delta}

\newcommand{\bD}{\boldsymbol D}
\newcommand{\bP}{\boldsymbol P}

\numberwithin{equation}{section}

\allowdisplaybreaks
\usepackage{xr}
\usepackage{url}
\begin{document}

\begin{frontmatter}
\title{Bayesian Linear Regression for Multivariate Responses Under Group Sparsity}
\runtitle{Bayesian variable selection for multivariate responses}

\begin{aug}
\author{\fnms{Bo} \snm{Ning}\thanksref{a}\ead[label=e1,mark]{bo.ning@yale.edu}},
\author{\fnms{Seonghyun} \snm{Jeong}\thanksref{b}\ead[label=e2,mark]{sjeong4@ncsu.edu}},
\and
\author{\fnms{Subhashis} \snm{Ghosal}\thanksref{b}\ead[label=e3,mark]{sghosal@ncsu.edu}\thanks{Research is partially supported by NSF grant number DMS-1510238.}}

\address[a]{Department of Statistics and Data Science, Yale University, 
24 Hillhouse Avenue, New Haven, CT 06511, USA}
\address[b]{Department of Statistics, North Carolina State University,
4276 SAS Hall, 2311 Stinson Drive, Raleigh,
NC 27695, USA
\printead{e1,e2,e3}}

\runauthor{Ning, Jeong, and Ghosal}

\affiliation{North Carolina State University}

\end{aug}

\begin{abstract}
We study frequentist properties of a Bayesian high-dimensional 
multivariate linear regression model with correlated responses.
The predictors are separated into many groups and the group structure is pre-determined.
Two features of the model are unique:
(i) group sparsity is imposed on the predictors.
(ii) the covariance matrix is unknown and its dimensions can also be high.
We choose a product of independent spike-and-slab priors 
on the regression coefficients 
and a new prior on the covariance matrix based on its eigendecomposition.
Each spike-and-slab prior is a mixture of a point mass at zero and a multivariate density involving a $\ell_{2,1}$-norm.
We first obtain the posterior contraction rate,
the bounds on the effective dimension of the model
with high posterior probabilities.
We then show that the multivariate regression coefficients
can be recovered under certain compatibility conditions.
Finally, we quantify the uncertainty for the regression coefficients
with frequentist validity through a Bernstein-von Mises type theorem.
The result leads to selection consistency for the Bayesian method.
We derive the posterior contraction rate using the general theory
by constructing a suitable test from the first principle using moment bounds for certain likelihood ratios.
This leads to posterior concentration around the truth with respect to 
the average R\'enyi divergence of order 1/2.
This technique of obtaining the required tests for posterior contraction rate could be useful in many other problems.

\end{abstract}

\begin{keyword}
\kwd{R\'enyi divergence}
\kwd{Bayesian variable selection}
\kwd{covariance matrix} 
\kwd{group sparsity} 
\kwd{multivariate linear regression}
\kwd{posterior contraction rate}
\kwd{spike-and-slab prior}
\end{keyword}

\end{frontmatter}

\section{Introduction}

Asymptotic behaviors of variable selection methods 
for linear regression were extensively studied
\citep{buhl11}.
However, theoretical studies on Bayesian variable selection methods
were limited to relatively simple settings
\citep{cast15, chae16, mart17, rock18, beli17b, song17}.
For example,
\citet{cast15} studied a sparse linear regression model in 
which the response variable is one-dimensional and the variance is known.
However, it is not straightforward to extend those results to multivariate linear regression with unknown covariance matrix (or even the univariate case with unknown variance).

Predictors can often be naturally clustered in groups, as in  the following examples.
\begin{itemize}
	\item[1.] {\it Cancer genomics study.} 
	The relationship between clinical phenotypes and DNA mutations is an important issue for biologists.
	DNA mutations are detected by DNA sequencing.
	Since these mutations are spaced linearly along the DNA sequence,
	it is often assumed that the adjacent DNA mutations on the chromosome have a similar genetic effect and should be grouped together \citep{li10}.
	
	\item[2.] {\it Multi-task learning.} 
	When information for multiple tasks is shared,
	solving tasks simultaneously is desirable   
	to improve learning efficiency and prediction accuracy.
	Relevant information is preserved across 
	different equations
	by grouping them together \citep{loun09}.
	
	\item[3.] {\it Causal inference in advertising.}
	When measuring the effectiveness of an advertising
	campaign running on stores,
	counterfactuals need to be constructed using the sales data at some control stores 
	chosen by a variable selection method
	\citep{ning18}.
	Stores within the same geographical region share the same demographic information, and so can be grouped together before selection. 
\end{itemize}

Driven by those applications,
new variable selection methods designed to select or not select variables as groups were developed by imposing {\it group-sparsity} on the regression coefficients as in the group-lasso \citep{ming06}.
This method replaces the $\ell_1$-norm in the penalty term of the lasso 
with the $\ell_{2,1}$-norm, which comprises of  
the $\ell_2$-norm put on the predictors within each group
and the $\ell_1$-norm is put across the groups.
Theoretical properties of the group-lasso were studied,
and its benefits over the lasso in the group selection problem were demonstrated 
\citep{nard08, loun09, loun11, huan10}.
Recently,
various Bayesian methods for selecting variables as groups 
were also proposed
\citep{li10, curt14, rock14, xu15, chen16, gree17, liqu17}.
However, their large-sample frequentist properties are largely unknown. 

In this paper, we study a Bayesian method for the multivariate linear regression model with two distinct features: group-sparsity imposed on the regression coefficients and an unknown covariance matrix.
To the best of our knowledge, even in a simpler setting without the group-sparsity structure, convergence and selection properties of methods for high-dimensional regression with a multivariate response having an unknown covariance matrix
have not been studied in either the frequentist or the Bayesian literature.
However, it is important to understand the theoretical properties of these methods because correlated responses arise in many applications.
For example, in the study of the causal effect of an advertising campaign, sales in different stores are usually spatially correlated \citep{ning18}.
Furthermore, when the dimension of the covariance matrix is large, 
it would affect the quality of the estimation of the regression coefficients.

When the covariance matrix is unknown and high-dimensional, 
standard techniques for posterior concentration rates \citep{cast15, mart17, beli17b} cannot be applied. 
Also, the general theory of posterior contraction under the average squared Hellinger distance \citep{ghos17} is not sufficient to obtain the rate in terms of the Euclidean metric on the regression parameter. In order to obtain that rate through the general theory, we shall construct certain required tests directly by controlling the moments of likelihood ratios with the parameter space broken up in small pieces. This leads to the posterior contraction rate with respect to the negative average log-affinity, which can be subsequently converted to the rate with respect to the Euclidean metric on the regression parameter.
The technique of controlling error probabilities by a moment bound on likelihood ratios appears to be new in the Bayesian literature and may be useful to study rates in other problems.

In this paper, we consider a multivariate linear regression model
\begin{equation}
\label{eqn-model}
Y_i = \sum_{j=1}^G X_{ij} \bbeta_j + \varepsilon_i,
\quad i = 1, \dots, n,
\end{equation}
where
$Y_i$ is a $1 \times d$ response variable, 
$i = 1, \dots, n$,
$X_{ij}$ is a $1 \times p_j$ predictor variable,
$j = 1, \dots, G$,
$\bbeta_j$ is a $p_j \times d$ matrix containing 
the regression coefficients,
and $\varepsilon_1, \dots, \varepsilon_n$ 
are independent identically distributed (i.i.d) as 
$\mathcal{N}(0, \bSigma)$ with
$\bSigma$ being a $d\times d$ unknown covariance matrix.
In other words, in the regression model, there are $G>1$ 
non-overlapping groups of predictor variables with the group structure being predetermined. 
We denote the groups which contain at least a non-zero coordinate as {\it non-zero groups} and the remaining groups as {\it zero groups}. 
The number of total groups $G$ is clearly bounded by $p$.
When $G = p$, it reduces to the setting that the sparsity is imposed on the individual coordinates.
Thus the results derived in our paper are also applicable to 
the ungrouped setting.

The above model can be rewritten in the vector form as 
\begin{equation*}
%\label{eqn-model-abbr}
Y_i = X_i \bbeta + \varepsilon_i,
\end{equation*}
where $\bbeta = (\bbeta_1', \dots, \bbeta_G')'$ 
is a $p \times d$ matrix, where $p = \sum_{j=1}^G p_j$,
and $X_i = (X_{i1}, \dots, X_{iG})$ is a $1 \times p$ vector.
%We will study the frequentist properties of 
%the posterior under the setting that 
The dimension $p$ can be very large and the dimension $d$ can be large as well, but to a lesser extent. 

To allow derivation of asymptotic properties of estimation and selection, certain conditions on the growth of $p$, $G$, $d$ and $p_1,\ldots,p_G$ need to be imposed. 
The dimension $p$ can grow at a rate faster than the sample size $n$, but we require that the total number of the coefficients in all non-zero groups together 
are less than $n$ in order. We further assume that
the number of coordinates in any single group must be of the order less than 
$n$, $G \geq n^c$, for some positive constant $c$, and $\log G$ grows slower than $n$.
Finally, to make the covariance matrix consistently estimable, we assume that for the dimension  $d$ of the covariance matrix, $d^2 \log n$ grows at a rate slower than $n$.

As for the priors,
we choose the product of $d$ independent 
spike-and-slab priors for $\bbeta$ and a prior on  
$\bSigma$ through its eigendecomposition. The latter seems to be a new addition to the literature.  
The spike-and-slab prior is a mixture of point mass for the 
zero coordinates
and a density for non-zero coordinates.
In the ungrouped setting, commonly used densities for 
non-zero coordinates are
a Laplace density \citep{cast15},
a Cauchy density \citep{cast18} and
a normal density with mean chosen by 
empirical Bayes methods \citep{mart17, beli17b}.
In this paper, we choose a density for the non-zero
coordinates involving the $\ell_{2,1}$-norm (see (\ref{eqn-prior-beta})), which corresponds to the penalty function of the group-lasso.
We derive an explicit expression for the normalizing constant of this density.

We shall use the following notations. We assume that $\mathcal{G}_1, \dots, \mathcal{G}_G$ 
are $G$ disjoint groups such that 
$\cup_{j=1}^G \mathcal{G}_j = \{1, \dots, p\}$.
Since these groups are given and will be kept the same throughout,
their notations will be dropped from subscription notations.
Each $p_j$ is the number of elements in $\mathcal{G}_j$.
Let $p_{\max} = \max\{p_j: 1 \leq j \leq G\}$.
For each $k = 1, \dots, d$,
let $S_k \subseteq \{1, \dots, G\}$
stand for the collection of indices of 
non-zero groups for the $k$th component
and $s_k = |S_k|$ stand its cardinality.
Let $S_{0,k}$ be the set consisting of 
the indices of the true non-zero groups.
Let
$S = \{S_1,\dots,S_d\}$ be the $d$-tuple of the model indices,
and define
$s = \sum_{k=1}^d s_k$,
$ 
p_{S_{k}} =  \sum_{j \in S_{k}} p_j,
$
and
$ 
p_{S} = \sum_{k=1}^d p_{S_{k}}.
$
Similar notations are used for the corresponding true values $S_{0,k}, s_{0,k}, S_0, s_0, p_{S_{0,k}}$ and 
$p_{S_0}$.
We also define $S_{\bbeta,k}$, $s_{\bbeta,k}$, $S_{\bbeta}$, and $s_{\bbeta}$ for an arbitrary $p\times d$ matrix $\bbeta$.

For a vector $A$, let $\|A\|_1$, $\|A\|_{2,1}$ and $\|A\|$ be the
$\ell_1$-, $\ell_{2,1}$- and $\ell_2$-norm of 
$A$, respectively, where
$\|A\|_{2,1} = \sum_{j=1}^G \|A_j\|$ with
$A_j$ being the subvector of $A$ consisting of $k \in \mathcal{G}_j$
coordinates.
For a matrix $\boldsymbol B$, 
let $\bmod{B}_{k}$ be the $k$th column of $\boldsymbol B$, by 
$\|\boldsymbol B\|_F =\sqrt{\mathrm{Tr}(\boldsymbol{B}^T \boldsymbol{B})}$ 
as the Frobenius norm, and $\lVert \boldsymbol{B} \rVert$ as the spectral norm. 
In particular, for an $n\times p$ matrix $\boldsymbol C$, we define the matrix norm $\lVert \boldsymbol{C} \rVert_\circ = \max\{\lVert \boldsymbol{C}_j \rVert: 1\le j\le G\}$, where $\boldsymbol{C}_j$  is the submatrix of $\boldsymbol{C}$ consisting of columns $C_k$ with $k \in \mathcal{G}_j$
coordinates.
For a $d \times d$ symmetric positive definite matrix $\boldsymbol D$,
let $\text{eig}_1(\boldsymbol D), \dots, \text{eig}_d(\boldsymbol D)$ denote the eigenvalues of $\boldsymbol C$ ordered from the smallest to the largest and 
$\det(\boldsymbol D)$ stands for the determinant of $\boldsymbol D$.
For a scalar $c$, 
we denote $|c|$ to be the absolute value of $c$.

Let $\rho(f, g) = -\log (\int f^{1/2} g^{1/2} d \nu)$ be the R\'enyi divergence of order 1/2 between densities $f$ and $g$
and 
$h^2(f, g) = \int (f^{1/2}- g^{1/2})^2 d\nu$ be their squared Hellinger distance.
The Kullback-Leibler divergence and the Kullback-Leibler variation between $f$ and $g$ are respectively given by 
$K(f, g) =\int f \log (f/g)$ and 
$V(f, g) = \int f (\log (f/g) - K(f, g))^2$.
The notation $\|\mu - \nu\|_{TV}$ denotes the total variation 
distance between two probability measures $\mu$ and $\nu$.

We let $N(\epsilon, \mathcal{F}, \rho)$ stand for the $\epsilon$-covering number of a set $\mathcal{F}$ with respect to a metric $\rho$, which is
the minimal number of $\epsilon$-balls in $\rho$-metric needed 
to cover the set $\mathcal{F}$.
Let $\bI_d$ stand for the $d$ dimensional identity matrix and
$\mathbbm{1}$ stand for the indicator function.

The symbols $\lesssim$ and $\gtrsim$ will be used to denote inequality up and down to a constant while
$a \asymp b$ stand for $C_1 a \leq b \leq C_2 a $ for two constants $C_1$ and $C_2$.
The notations $a \ll b$ and $a \vee b$ stand for $a/b \rightarrow 0$ and
$\max\{a, b\}$ respectively.
The symbol $\delta_0(\cdot)$ stands for the probability measure with all its mass at $0$.

The remainder of the paper is organized as follows.
Section \ref{sec-priors} describes the priors, 
along with the necessary assumptions.
Section \ref{sec-main-results} provides the main results.
Section \ref{sec-computation} discusses algorithms for computation.
The proofs of two main theorems are given in Section \ref{sec-proof}.
The supplementary material gives an auxiliary result and presents analogous but slightly weaker results on posterior contraction and selection using a conjugate inverse-Wishart prior on the covariance matrix.

\section{Prior specifications}
\label{sec-priors}

In this section, we introduce the priors used in this paper.
We let $\bbeta$ and $\bSigma$ be independently distributed in the prior. The prior for $\bbeta$ is mixed over several dimensions and each component of the prior density depends on the $\ell_{2,1}$-norm of $\bbeta$, while a spike-and-slab prior is put on the group dimension. We put a prior on the covariance matrix through its eigendecomposition $\bSigma = \bP \bD \bP'$, with  independent inverse Gaussian priors for each diagonal entry of $\bD$ and the uniform prior for $\bP$ on the group of orthogonal matrices.

\subsection{Prior for regression coefficients}

We denote the $k$th column of $\bbeta$ by $\beta_k$ and let 
the notations
$\beta_{k,S_k}$ and $\beta_{k,S_k^c}$ stand for the  collection of regression coefficients in the $k$th column of the non-zero groups and the zero groups respectively.
A spike-and-slab prior is constructed as follows.
First, the dimension $s$ is chosen from a prior $\pi$ on the set $\{0, 1, \dots, Gd\}$.
Next, a $d$-tuple $S$ of subsets is randomly chosen from the set $\{1, \dots, G\}^d$ such that $\sum_{k=1}^d s_k=s$.
Finally, for each $k$, a vector 
$\beta_{k,S_k}$ 
is independently chosen from a
probability density $g_{S_k}$ on $\mathbb{R}^{p_{S_k}}$ given by 
(\ref{eqn-gbeta}), and the remaining coordinates $\beta_{k,S_k^c}$ set to 0.
To summarize,
the prior for $\bbeta$ is
\begin{align}
\label{eqn-prior-beta}
(S, \bbeta)  \rightarrow \pi(s) \frac{1}{{Gd \choose s}} 
\prod_{k=1}^d g_{S_k}(\beta_{k,S_k}) 
\delta_0(\beta_{k,S_k^c}),
\end{align}
where the density $\pi(s)$ is the prior for the dimension $s$.

\begin{assumption}[Prior on dimension]
	\label{assump-prior-dimension}
	For some constants $A_1, A_2, A_3, A_4>0$,  
	\begin{align*}
	%\label{eqn-assump-dimension}
	\frac{A_1}{(G\vee n^{p_{\max}})^{A_3}} 
	\leq 
	\frac{\pi(s)}{\pi(s - 1)}
	\leq
	\frac{A_2}{(G\vee n^{p_{\max}})^{A_4}}  ,\quad s = 1, \dots, Gd.
	\end{align*}
\end{assumption}
If sparsity is imposed at the individual level, i.e. $p_{\max}=1$, then the assumption is identical to the one given in \citet{cast15}. Prior distributions satisfying the assumption can easily be constructed.
For example, the complexity prior given by \citet{cast15} satisfies the above assumption if $p_{\max}=1$, and it can also be easily modified to consider the case when $p_{\max}>1$.

When sparsity is at the individual level,
the Laplace density \citep{cast15} or the Cauchy density \citep{cast18} 
is generally chosen for $g$, since the normal density has a too sharp tail that overshrinks the non-zero coefficients, although some empirical Bayes modifications of the mean can overcome the issue
\citep[see][]{mart17, beli17b}.
However, in our setting, 
as sparsity is imposed at the group level,  like the group lasso, 
we consider the following density using the $\ell_{2,1}$-norm:
\begin{align}
\label{eqn-gbeta}
g_{S_k}(\beta_{k,S_k}) = \Bigg(
\prod_{j \in S_k} 
\Big(\frac{\lambda_k}{a_j}\Big)^{p_j}
\Bigg)
\exp\big(-\lambda_k \|\beta_{k,S_k}\|_{2,1}\big),
\end{align}
where
$
a_j =
\sqrt{\pi}
({\Gamma(p_j+1)}
/{\Gamma(p_j/2 + 1)})^{1/p_j} \geq 2$
(see Lemma~\ref{lemma-apx-const-a} in the supplementary materials).
%\ref{lemma-apx-const-a}
This density has its tail lighter than the corresponding
Laplace density.
From Stirling's approximation,
it follows that
$a_j = O(p_j^{1/2})$. A relevant elliptical prior distribution is considered in \citet{gao17}.

A prior of this type involving the $\ell_{2,1}$-norm was also used in the Bayesian literature in group-sparsity problems \citep{xu15}, but an explicit expression of the normalizing constant was not obtained. 
Since the normalizing constant depends on the dimension, its value will play a role in the posterior contraction rate.

The tuning parameter $\lambda_k$ in the prior needs to be bounded both from above and below, specified in Assumption~\ref{assump-lambda} below. A value too large will shrink the non-zero coordinates too much towards to 0.
A value too small will be unable to prevent many false signals appearing in the model, which can make the posterior to contract slower. 

\begin{assumption}
	\label{assump-lambda} 
	For some constants $B_1, B_2, B_3>0$ and each $k = 1, \dots, d$, 
	$\underline \lambda \leq \lambda_k \leq \overline\lambda$,
	where
	\begin{align}
	\label{eqn-assump-2}
	\underline{\lambda} =\frac{\lVert \bX\rVert_\circ}{B_1(G^{1/p_{\max}}\vee n)^{B_2}}
	\quad
	\overline{\lambda}
	= B_3\lVert \bX\rVert_\circ \sqrt{\log G\vee p_{\max}\log n }.
	\end{align}
	
\end{assumption}

The constants $B_1$, $B_2$, $B_3$ can be chosen large enough so that the range can be sufficiently wide. In particular, if $p_{\max}=1$, this above reduces to the one in \citet{cast15}.

Assumption~\ref{assump-lambda} will be coupled with Assumption~\ref{assump-true-parameters} in Section~\ref{sec-post-rate} on the true parameters. A particularly interesting case is that every $\lambda_k$ is set to the lower bound $\underline \lambda$ for every $k$. Then the bound requirement on the true signal will be rather mild.

\subsection{Prior for the covariance matrix}
\label{prior-covariance}

For a prior on the covariance matrix $\bSigma$, we use its eigendecomposition $\bP \bD \bP'$. We put an inverse Gaussian prior independently on each eigenvalue of $\bSigma$, or equivalently, on each diagonal entry of $\bD$. This prior is chosen because of its exponentially decaying tail on both sides. The orthogonal matrix $\bP$ is given the Haar measure on the compact Lie group of $d\times d$ orthogonal matrices, which is a Riemannian manifold of  dimension $d(d-1)/2$ embedded in $\mathbb{R}^{d\times d}$.

We found that the naturally conjugate inverse Wishart prior on $\bSigma$ may induce a suboptimal posterior contraction rate due to its weaker tail property when $d$ increases to infinity. Nevertheless, because of the practical importance of this prior, we present the contraction rate for this prior in the supplementary material.  
When $d$ is fixed, the rate is the same as in the main theorem in this paper using the above stated prior on $\bSigma$.   
When additional structure like sparsity are assumed on large covariance or precision (inverse covariance) matrices, prior distributions can be assigned by respecting such structure \citep{bane14,bane15,pati2014posterior}. In such a situation, an improved rate may be possible; see the remark at the end
of Section~\ref{sec-post-rate}.
Other significant priors used in the literature, 
such as reference priors \citep{yang94, sun07}, are harder to handle because the general theory of posterior contraction does not apply to these improper priors, and moreover, tail bounds for the corresponding eigenvalue distribution need to be available.

\section{Main results}
\label{sec-main-results}

\subsection{Posterior contraction rate}
\label{sec-post-rate}

We study the posterior contraction rate 
for the model
and the priors given in Section \ref{sec-priors}.
We denote $\bbeta_0$ and $\bSigma_0$ as the true values
of $\bbeta$ and $\bSigma$, respectively.
Recall the notations
$s_{0,k} = |S_{0,k}|$,
$S_0 = \{S_{0,1},\dots,S_{0,d}\}$,
and $s_0 = \sum_{k=1}^d s_{0,k}$.

The general theory of posterior contraction for independent non-identically distributed observations \citep[see Theorem 8.23 of][]{ghos17} is often used to derive a posterior contraction rate. The general theory characterizes the contraction rate in terms of the average squared Hellinger distance by default, unless an additional testing property in the model is established. However, closeness in terms of the average squared Hellinger distance between multivariate normal densities with varying mean and an unknown covariance does not necessarily imply that the mean parameters in the two densities are also close on average in terms of the Euclidean distance. To alleviate the problem, we work directly with the average R\'enyi divergence of order 1/2, which is still very tractable in the multivariate normal setting, and gives rise to closeness in terms of the desirable Euclidean distance. To this end, we directly construct a suitable test using the likelihood ratio for the null against some representative points in the alternative described by the complement of a R\'enyi ball around the null intersected with a sieve, and then showing that such a test also works well for testing the null value against a neighborhood of the representative point, by controlling the moments of the likelihood ratio of the representative point and the points in the neighborhood. Finally, by controlling the number of pieces needed to cover the sieve, we construct a single test with required control over the error probabilities for testing the null value against the whole of the alternative intersected with the sieve, which can then be used in the general theory of posterior contraction. 

The general theory also requires lower bounds for prior concentration near the true parameter value, which is possible provided that we require the true values of $\bbeta_0$ and $\bSigma_0$ to be restricted into certain regions (see Assumption~\ref{assump-true-parameters} below). 
This is unlike \citet{cast15}, who obtained results uniformly over the whole space as their case  
(univariate with known variance and Laplace prior)
allows explicit expressions for direct treatment.

\begin{assumption}
	\label{assump-true-parameters}
	The true values satisfy $\bbeta_0 \in \mathcal{B}_0$ 
	and $\bSigma_0 \in \mathcal{H}_0$, for
	\begin{eqnarray}
	\mathcal{B}_0 = 
	\bigg\{\bbeta:
	\sum_{k=1}^d
	\|\beta_{k}\|_{2,1} 
	\leq \overline{\beta}\bigg\}, 
	\quad 
	\mathcal{H}_0 = 
	\{\bSigma: 
	b_1 \bI_d \leq \bSigma \leq b_2 \bI_d\},
	\label{eqn-true-parameters}
	\end{eqnarray}
	where $b_1,b_2>0$ are fixed values and
	$\overline{\beta} = 
	{ s_0
		(\log G 
		\vee 
		{p_{\max} \log n })}/
	{\max\{ \lambda_k: 1\leq k\leq d\} }
	$.
\end{assumption}

The largest value of $\overline{\beta}$ is obtained by taking $\lambda_k = \underline{\lambda}$ for all $k$.
In this case, the upper bound becomes $\overline{\beta} = {B_1 s_0 (\log G\vee p_{\max}\log n)(G^{1/p_{\max}}\vee n)^{B_2}}/{\lVert \bX\rVert_\circ}$, which is a very mild restriction if $B_2$ is chosen large enough. 

\begin{theorem}
	\label{thm-3.1}
	For the model (\ref{eqn-model}) and the priors
	given in Section \ref{sec-priors}, we have that for a sufficiently large $M_1>0$,
	\begin{align}
	\sup_{\bbeta_0 \in \mathcal{B}_0,
		\bSigma_0 \in \mathcal{H}_0}
	\mathbb{E}_0  \Pi\Big(
	\bbeta: 
	\|\bX(\bbeta - \bbeta_0)\|_F^2
	\geq M_1 n\epsilon_n^2 \Big| Y_1, \dots, Y_n
	\Big) &\rightarrow 0,
	\label{eqn-post-beta}\\
	\sup_{\bbeta_0 \in \mathcal{B}_0,
		\bSigma_0 \in \mathcal{H}_0}
	\mathbb{E}_0  \Pi\Big(
	\bSigma:
	\|\bSigma - \bSigma_0\|_F^2
	\geq M_1 \epsilon_n^2 \Big| Y_1, \dots, Y_n
	\Big) &\rightarrow 0,
	\label{eqn-post-sigma}
	\end{align}
	where
	\begin{align}
	\label{eqn-post-rate}
	\epsilon_n = \max \left \{
	\sqrt{\frac{s_{0}  \log G}{n}} 
	,
	\sqrt{\frac{s_0 p_{\max}  \log n}{n}}
	,
	\sqrt{\frac{d^2\log n}{n}} \right\}\to 0.
	\end{align}
\end{theorem}

\begin{remark}
	\rm
	Unlike in the classical approach where variable selection is regulated by a penalty function that corresponds to a prior density on the regression coefficients, in the Bayesian approach, sparsity is imposed by the prior on the dimension. The prior density on the regression coefficients still plays a significant role in controlling the prior concentration and the tail behavior, but to a lesser extent. Thus,
	instead of using the prior given in (\ref{eqn-gbeta}),
	one can also choose a Laplace density for the coordinates in 
	the non-zero groups. 
	Then the $\ell_{2,1}$-norm of $\beta_{0,k}$,
	$\|\beta_{0,k}\|_{2,1}$, in the set $\mathcal{B}_0$ should be replaced by $\|\beta_{0,k}\|_1$.
	Clearly,
	$\|\beta_{0,k}\|_{2,1} \leq \|\beta_{0,k}\|_1$,
	and hence in the latter case, 
	the set $\mathcal{B}_0$ will be smaller.
\end{remark}

\begin{remark}
	\rm
	When $G = p$, and hence $p_{\max}=1$, 
	the posterior contraction rate simplifies to 
	$\epsilon_n = \max \{
	\sqrt{{(s_0 \log p)}/{n}}, 
	\sqrt{{(d^2\log n)}/{n}}\}$.
	The first term in the rate is the same as the rate obtained 
	when the sparsity is imposed at the individual level,
	such as in \citet{buhl11} and \citet{cast15}.
	When $G \ll p$,
	the same rate can be obtained if
	$\displaystyle
	p_{\max}\log n \lesssim \log G$.
	%(i.e., when the number of coordinates in each group 
	% takes a fixed number) and $d$ is sufficiently slowly growing.
\end{remark}

The first term of the rate in Theorem \ref{thm-3.1} coincides with the rate obtained for a group-lasso estimator of the multi-task learning problem studied by \citet{loun11}. Their setup is not directly comparable with ours but their analogous rate coincides with ours up to a logarithmic factor and they showed its optimality in a minimax sense. Uder the setting $d=1$, the rate obtained in \citet{huan10} is  $(p_{S_0} + s_0 \log G)/n$, which is only slightly faster than our rate, and will coincide with ours up to the logarithmic factor whenever  $p_{S_0} \asymp s_0 p_{\max}$. This can often happen provided that the non-zero groups are not consisting of a few large and the rest small groups. 

If there is additional lower-dimensional structure in the orthogonal matrix $\bm{P}$, the last term in \eqref{eqn-post-rate} may be improved, because in a lower-dimensional manifold, the prior concentration rate will be higher and the entropy estimates will be lower. The simplest such structure is the trivial situation $\bP=\bI$, which leads to diagonal covariance matrix and the  reduction of $d^2$ to $d$. More generally, a block-diagonal structure with $L$ non-overlapping blocks of size $d_1,\ldots,d_L$, $\sum_{l=1}^L d_l=d$, will reduce $d^2$ to $\sum_{l=1}^L d_l^2$.

From Theorem \ref{thm-3.1}, the posterior contraction rate slows down significantly if the dimension of the covariance 
is too high, but a better rate may be possible if a lower dimensional structures is present in the covariance of the precision matrix. For instance, if the responses are independent across components,
then the model (\ref{eqn-model}) can be written as $d$ independent model
with each one is
\begin{align*}
{\sigma_k}^{-1} Y_{ik} = {\sigma_k}^{-1} X_{i} \beta_k
+ \varepsilon_{ik}, \quad \varepsilon_{ik} \sim \mathcal{N}(0, 1).
\end{align*}
Then one can estimate the parameters in the $d$ models separately.
The posterior concentration rate for each corresponding posterior becomes
$\epsilon_{n} = (\sum_{k=1}^d \epsilon_{n,k}^2 )^{1/2}$,
where
$\epsilon_{n,k} = 
\max \{
\sqrt{({s_{0,k} \log G})/{n}}
, \sqrt{(
	s_{0,k}
	p_{\max} \log n)/n}
\}$  is the individual rates for the $k$th component, $k =1, \dots, d$. 

\subsection{Dimensionality and recovery}
\label{sec-recovery}

In this section, 
we show dimensionality control and recovery properties of the
the marginal posterior of $\bbeta$.

\begin{lemma}[Dimension]
	\label{lemma-dim}
	For the model (\ref{eqn-model}) and the priors
	given in Section \ref{sec-priors}, we have that for a sufficiently large number $M_2 >0$,
	\begin{equation*}
	%\label{equ-dim}
	\sup_{\bbeta_0 \in \mathcal{B}_0,
		\bSigma_0 \in \mathcal{H}_0}
	\mathbb{E}_0
	\Pi \Big(\bbeta: s_{\bbeta}
	\geq
	M_2 s^\star
	\Big | Y_1, \dots, Y_n
	\Big) \rightarrow 0,
	\end{equation*}
	where $ s^\star = s_0\vee \{ {d^2\log n}/({\log G\vee p_{\max}\log n})\}$.
\end{lemma}

From Lemma \ref{lemma-dim},
$s^\star > s_0$ if
${d^2 \log n}  \gg {s_0(\log G\vee  p_{\max}\log n)}$.
This means that the support of the posterior can substantially 
overshoot the true dimension $s_0$.
In the next corollary, we show that even when $s^\star > s_0$, the posterior is still able to recover $\bbeta_0$ in terms of the distance to the truth. 

\begin{corollary}[Recovery]
	\label{cor-recovery}
	For the model (\ref{eqn-model}) and the priors
	given in Section \ref{sec-priors}, we have that for a
	sufficiently large constant $M_3 > 0$,
	\begin{align}
	\sup_{\bbeta_0 \in \mathcal{B}_0,
		\bSigma_0 \in \mathcal{H}_0}
	\mathbb{E}_0
	\Pi
	\left(
	\|\bbeta - \bbeta_{0}\|_F^2
	\geq 
	\frac{M_3 n\epsilon_n^2 }
	{
		\lVert \bX\rVert_\circ^2
		\phi_{\ell_2}^2(s_0+M_2s^\star)
	}
	\Bigg| Y_1, \dots, Y_n
	\right) 
	\rightarrow 0,
	\label{eqn-recovery-l2}
	\end{align}
	where $\phi_{\ell_2}^2$ is the restricted eigenvalue
	(see Definition \ref{def-res} below).
\end{corollary}

\begin{definition}[Restricted eigenvalue]
	\label{def-res}
	The smallest scaled singular value of dimension $\tilde{s}$ is defined as
	\begin{equation}
	\label{assump-res}
	\phi_{\ell_2}^2(\tilde{s}) = 
	\inf\Bigg\{
	\frac{ 
		\|\bX\bbeta\|_F^2
	}{
		\lVert \bX\rVert_\circ^2
		\|\bbeta\|_F^2}
	, \ 0 \leq s_{\bbeta} \leq \tilde{s}
	\Bigg\}.
	\end{equation}
\end{definition}   

As $p \gg n$, 
the smallest eigenvalue of the design matrix must be 0.
The restricted eigenvalue condition keeps the smallest eigenvalue
for the sub-matrix of the design matrix, corresponding  
to the coefficients within non-zero groups, bounded away from 0. 

The results in terms of other norms for the difference between $\bbeta$ and $\bbeta_0$ can be also derived by assuming different assumptions
on the smallest eigenvalue for the sub-matrix of the design matrix.
For example, by using the uniform compatibility condition (in Definition \ref{def-comp-l21} below),
we can conclude that for a sufficiently large number $M_4 > 0$,
\begin{align}
\label{eqn-recover-l2l1}
\sup_{\bbeta_0 \in \mathcal{B}_0,
	\bSigma_0 \in \mathcal{H}_0}
\mathbb{E}_0
\Pi
\left( \bigg(
\sum_{k=1}^d
\|\beta_k - \beta_{0,k}\|_{2,1}\bigg)^2
\geq 
\frac{M_4 s^\star n\epsilon_n^2  }
{
	\lVert \bX\rVert_\circ^2
	\phi_{\ell_{2,1}}^2(s_0+M_2s^\star)
}
\Bigg| Y_1, \dots, Y_n
\right) 
\rightarrow 0.
\end{align}
We omit the proof since it 
is almost identical to that of Corollary \ref{cor-recovery}.

\begin{definition}[Uniform compatibility, $\ell_{2,1}$-norm]
	\label{def-comp-l21}
	The $\ell_{2,1}$-compatibility number in vectors of dimension $\tilde{s}$ is defined as
	\begin{equation*}
	%\label{assump-comp-l21}
	\phi_{\ell_{2,1}}^2(\tilde{s}) = 
	\inf\Bigg\{
	\frac{
		s_{\bbeta}
		\|\bX \bbeta\|_F^2
	}{
		\lVert \bX\rVert_\circ^2
		(\sum_{k=1}^d \|\beta_k\|_{2,1})^2
	}, \
	0 \leq s_{\bbeta} \leq \tilde{s}
	\Bigg\}.
	\end{equation*}
\end{definition}
By the Cauchy-Schwarz inequality,
$\sqrt{s_{\bbeta}} \|\bbeta\| _F
\geq \sum_{k=1}^d{\|\beta_k\|_{2,1}}$,
and it follows that
$\phi_{\ell_2}(\tilde{s}) \leq \phi_{\ell_{2,1}}(\tilde{s})$ for any 
$\tilde{s} \ll Gd$.

\subsection{Distributional approximation}
\label{sec-bvm}

To establish selection consistency, \citet{cast15} devised a key technique through a distributional approximation for the posterior distribution. As in a Bernstein-von Mises (BvM) theorem, the posterior distribution of the regression parameter is approximated by a relatively simpler distribution, but unlike in a traditional BvM theorem for increasing dimensional parameters \citep{ghosal1999asymptotic,ghosal2000asymptotic,bontemps2011bernstein} or low-dimensional functionals \citep{deJong13,gao16}, the approximating distribution is a mixture of multivariate normal instead of a single one.  

To derive an appropriate distributional approximation, we rewrite the model (\ref{eqn-model}) as 
\begin{align*}
%\label{eqn-model-2}
Y_i = \vc(\bbeta)\tilde\bX_i + \varepsilon_i, 
\quad 
i = 1, \dots, n,
\end{align*}
where $\vc(\bbeta)$ is obtained 
by stacking all the columns of $\bbeta$ 
into a $pd$-dimensional row vector,
$\tilde\bX_i = \bI_d \otimes X_i'$
is a $pd \times d$ block-diagonal matrix. The log-likelihood function is given by
\begin{equation}
\begin{split}
\ell_n( \bbeta , \bSigma)
%& =     \sum_{i=1}^n 
% \log f\big(Y_i|\vc(\bbeta)\bX_{i}, \bSigma\big)\\
& =
-\frac{nd}{2}\log(2\pi)
-\frac{n}{2}\log\big(\det(\bSigma)\big)
-\frac{1}{2}\sum_{i=1}^n \lVert\bSigma^{-1/2}
\big(Y_i - \vc(\bbeta) \tilde\bX_{i} \big)'\rVert^2.
\label{eqn-likelihood}
\end{split}
\end{equation}
For any measurable subset $\mathcal{B}$ of $\mathbb{R}^{p\times d}$, the marginal posterior distribution of $\bbeta$ is 
\begin{align}
\label{eqn-post-mixture}
\Pi(\bbeta\in\mathcal{B}|Y_1, \dots, Y_n)
=
\frac{
	\int \int_{\mathcal{B}} 
	\exp\big(\ell_n( \bbeta , \bSigma) -
	\ell_n( \bbeta_0, \bSigma_0) \big) 
	d\Pi ( \bbeta)d\Pi(\bSigma) 
}
{
	\int \int 
	\exp\big(\ell( \bbeta , \bSigma) -
	\ell( \bbeta_0, \bSigma_0) \big) 
	d\Pi ( \bbeta)d\Pi(\bSigma) },
\end{align}
with 
\begin{align*}
d\Pi ( \bbeta )
= 
\sum_{S:s\le Gd}
\frac{\pi(s)}{{Gd \choose s}}\prod_{k=1}^d \left\{
\Bigg( \prod_{j \in S_k}\Big(\frac{\lambda_k}{a_j}\Big)^{p_j} \Bigg)
\exp (- \lambda_k\|\beta_{k,S_k}\|_{2,1})
d\beta_{k,S_k} \otimes \delta_{S_{k}^c}\right\}.
\end{align*}

In the next theorem, we shall show that under certain conditions,
the posterior probability $\Pi(\bbeta\in\mathcal{B}|Y_1, \dots, Y_n)$ 
can be approximated by
\begin{align*}
\Pi^\infty(\bbeta \in \mathcal{B}|Y_1, \dots, Y_n)
=
\frac{
	\int_{\mathcal{B}} 
	\exp\{\ell_n (\bbeta, \bSigma_0 )
	- \ell_n (\bbeta_0, \bSigma_0)
	\} dU(\bbeta)
}
{
	\int
	\exp\{\ell_n (\bbeta, \bSigma_0 )
	- \ell_n (\bbeta_0, \bSigma_0)
	\} dU(\bbeta)
},
\end{align*}
where
\begin{align}
dU ( \bbeta) = 
\sum_{S : s\le M_2 s^\star}
\frac{\pi(s)}{{Gd \choose s}}
\prod_{k=1}^d \left\{\Bigg(\prod_{j\in S_k}\left(\frac{\lambda_k}{a_j}\right)^{p_j}\Bigg) d \beta_{k,S_k}\otimes\delta_{S_k^c}\right\}.
\label{eqn:appdu}
\end{align}
This means that $\ell_n( \bbeta , \bSigma)$ can be replaced by $\ell_n( \bbeta , \bSigma_0)$ with the true $\bSigma_0$ and the impact of the $\ell_{2,1}$-term in the prior density vanishes. Let
$\tilde\bX_{i,S}$ be the submatrix of $\tilde\bX_i$ chosen by $S$, with its dimension $p_S \times d$.
If $p_S\le n$ for a given $S$,
the maximum likelihood estimator (MLE) for $ \beta_S^\star = (\beta_{1,S_1}',\dots,\beta_{d,S_d}')'$ given the true covariance matrix $\bSigma_0$
is unique.
% Define $Y_i^S = Y_i - \vc(\bbeta_{S^c}) \bX_{i,S^c}$,
We denote the MLE and the information matrix as 
\begin{align*}
\hat\beta_S^\star & = 
\left( \sum_{i=1}^n \tilde\bX_{i,S} \bSigma_0^{-1} \tilde\bX_{i,S}' \right)^{-1}
\left(\sum_{i=1}^n \tilde\bX_{i,S} \bSigma_0^{-1} Y_i' \right),\quad
\hat{\mathbb{I}}_S = \frac{1}{n}\sum_{i=1}^n 
\tilde\bX_{i,S}\bSigma_0^{-1}\tilde\bX_{i,S}'.
\end{align*}
Then we can also write
\begin{align}
\label{eqn-post-mixture-infty}
\Pi^{\infty}(\bbeta\in\cdot|Y_1, \dots, Y_n)
& \propto
\sum_{S:s\le M_2 s^\star} 
w_S^\infty
\mathcal{N}\left(\hat\beta_S^\star , 
n^{-1}\hat{\mathbb{I}}_{S}^{-1}\right) 
\otimes \delta_{S^c},
\end{align}
where
\begin{align*}
%\label{eqn-post-mixture-infty-weights}
w_S^\infty
& \propto
\frac{ \pi(s)}{{Gd \choose s}}
\Bigg(\prod_{k=1}^d \prod_{j \in S_k} \bigg(\frac{\lambda_k\sqrt{2\pi}}{a_j}\bigg)^{p_j}\Bigg){\det \big(n \hat{\mathbb{I}}_{S}\big)^{-1/2}}
\exp\left\{\frac{1}{2}\sum_{i=1}^n
\|  \bSigma_0^{-1/2}  \bX_{i,S}' \hat{\beta}_S^\star\|^2\right\},
\end{align*}
with $\sum_S w_S^\infty = 1$.

Before we formally state the theorem, we recall the notion of the {\it small $\lambda$ regime} \citep[see][]{cast15}. Clearly, bounded $\lambda$-values belong to the small $\lambda$ regime.  
In our setting, we say $\lambda_k$ belongs to the {\it small $\lambda$ regime} if $\max\{
\lambda_k \epsilon_n\sqrt{s^\star n}/\lVert \bX\rVert_\circ: 1 \leq k \leq d\} \rightarrow 0$.
In this regime, the impact of the $\ell_{2,1}$-penalty vanishes, and hence
the MLE $\hat\beta_S^\star$ is asymptotically
unbiased and does not depend on the choice 
of different values of $\lambda_k$.
When choosing the value of $\lambda_k$ outside the small $\lambda$ regime,
this MLE is no longer asymptotically unbiased (see Theorem 11 of the supplementary material of \citet{cast15}). In order to make the remainder of the approximation tend to zero, we also assume that $\epsilon_n^2\sqrt{s^\star n}\big(\sqrt{s^\star n\epsilon_n^2 }\vee\sqrt{p_{\max} d^3  \log G}\big)\rightarrow 0$ for  $\ell_n( \bbeta , \bSigma)$ to be replaced by $\ell_n( \bbeta , \bSigma_0)$.

\begin{theorem}[Distributional approximation]
	\label{thm-bvm}
	For the model (\ref{eqn-model}), the priors
	given in Section \ref{sec-priors} with $\lambda$ in the small $\lambda$ regime,  and the sequence 
	\begin{align*}
	\delta_n(s_0)=\epsilon_n\sqrt{s^\star n}\max\left({\max\{
		\lambda_k: 1 \leq k \leq d\}  }/{\lVert \bX\rVert_\circ},\epsilon_n^2\sqrt{s^\star n },\epsilon_n\sqrt{p_{\max} d^3  \log G}\right),
	\end{align*}
	we have that any positive sequence $\eta_n\to 0$ and some positive constant $c>0$,
	\begin{align*}
	%\label{eqn-thm4.2}
	\sup_{\substack{
			\bbeta_0 \in  \{\mathcal{B}_0: \delta_n(s_0)<\eta_n,\\ \phi_{\ell_{2,1}}(s_0+M_2s^\star) > c \} ,\,
			\bSigma_0 \in \mathcal{H}_0
	}}
	\mathbb{E}_0
	\|\Pi(\bbeta\in\cdot|Y_1, \dots, Y_n) 
	- \Pi^{\infty}(\bbeta\in\cdot|Y_1, \dots, Y_n)\|_{TV}
	\rightarrow 0.
	\end{align*}
\end{theorem}

\subsection{Selection}
\label{sec-selection}

In this section, we establish selection consistency using Bernstein-von Mises theorem of the previous section. We assume the dimension of the covariance 
and the coordinates in the non-zero groups are sufficiently small.
We also assume the smallest signal cannot be too small,
which is
\begin{align}
\label{eqn-beta-min}
\tilde{\mathcal{B}} = 
\left\{\bbeta:
\min\{
\|\beta_{jk}\|^2: j \in S_{0,k},
k=1,\ldots, d\} \geq 
\frac{M_3 n\epsilon_n^2 }
{\lVert \bX\rVert_\circ^2
	\phi_{\ell_2}^2(s_0+M_2s^\star)}
\right\}.
\end{align}
This condition can be viewed as the {\it Beta-min condition} under the group sparsity setting.
The lower bound displayed in the condition is derived from (\ref{eqn-recovery-l2}).
Unlike the {\it Beta-min condition} in \citet{cast15}
which requires each individual coordinate is
bounded away from 0, our condition allows 
zero coordinates to be included in a non-zero group. 

The Beta-min condition is not vacuous, in that the lower bound in (\ref{eqn-true-parameters}) is smaller than the upper bound in (\ref{eqn-beta-min}). To see this, note that under the small $\lambda$ regime, $(\max_k \lambda_k)^{-1} \gg \sqrt{s^\star n\epsilon_n^2}/\|\bX\|_\circ$. Therefore, $\bar{\beta}\gg \sqrt{n \epsilon_n^2}/\|\bX\|_\circ$, and  the right side coincides with the lower bound up to a constant, establishing the claim.

We now complete this section by stating the following theorem.

\begin{theorem}[Selection consistency]
	\label{thm-selection}
	For the model (\ref{eqn-model}), the priors
	given in Section \ref{sec-priors}, some positive constant $c>0$, and some sequences $\eta_n\rightarrow 0$ and $s_n\le G^a$ with $a<A_4-3/2$,
	we have that
	\begin{align*}
	\sup_{\substack{
			\bbeta_0 \in  \{\mathcal{B}_0\cap  \tilde{\mathcal{B}}: s_0\le s_n, \delta_n(s_0)\le \eta_n,\\ \phi_{\ell_{2,1}}(s_0+M_2s^\star) > c \} ,\,
			\bSigma_0 \in \mathcal{H}_0
	}}
	\mathbb{E}_0\Pi(\bbeta: S_{\bbeta} = S_0
	|Y_1, \dots, Y_n) 
	\rightarrow 1.
	\end{align*}
\end{theorem}

Under the conditions in Theorem \ref{thm-selection},
the marginal posterior distribution of $\bbeta$ in non-zero groups can be further approximated by a multivariate normal distribution
with mean $\vc(\hat{\beta}_{S_0}^\star)$
and the covariance matrix 
$\hat{\mathbbm{I}}_{S_0}^{-1}
= n \big(\sum_{i=1}^n 
\bX_{i,S_0}\bSigma_0^{-1}\bX_{i,S_0}')^{-1}$. 
Therefore, credible sets for $\bbeta$ can be obtained directly from the approximating multivariate normal density.
It may be noted that under the setting of the theorem, the lower bound in the Beta-min condition goes to zero, implying that the condition becomes milder with increasing sample size.

\section{Computational algorithms}
\label{sec-computation}

Various sampling-based computation algorithms have been developed to compute the posterior distribution in the sparse linear regression model with a spike-and-slab prior under the setting that the covariance matrix $\bSigma = \sigma^2 \bI_d$ and sparsity is imposed on individual coefficients.
A summary of those algorithms is provided in Section 5 of \citet{cast15}.
Recently, \citet{xu15} developed an MCMC algorithm using a spike-and-slab prior for group variable selection. They placed a beta-binomial prior on the dimension and a prior involves $\ell_{2,1}$ norm, similar to ours, on the regression coefficients.

Since priors used in this paper are new,
we outline an MCMC algorithm to compute the posterior distribution.
For each iteration of the algorithm, one can start with sampling $S$ from the marginal posterior distribution with $\bbeta$ integrated out. Next, conditioning on the current $S$, draw $\bbeta$ from the corresponding conditional posterior distribution. Since the prior for $\bbeta$ is not a conjugate prior, the Metropolis-Hasting algorithm can be used with the proposal density chosen as a multivariate normal distribution centered at its current value.
Last, sample $\bSigma$ through sampling $\boldsymbol P$ and $\boldsymbol D$, and then calculating $\boldsymbol P \boldsymbol D \boldsymbol P'$.
To sample the diagonal elements of $\boldsymbol D$, one can convert them to log scale and then for each element, choose the proposal density as a normal distribution centered at its current value in $\log$ scale. 
To sample $\boldsymbol P$, one can draw a new value $\boldsymbol P^\star$ uniformly from the group of orthogonal matrices.
Then the acceptance ratio equals to the likelihood ratio.
When $d$ is large, in order to increase acceptance rate of the Metropolis-Hasting algorithm,
one can restrict the proposal density to local moves through multiplying by a random orthogonal matrix with some $\epsilon$ of the identity matrix.
If the conjugate inverse Wishart prior is used instead, then the conditional posterior distribution of $\bSigma$ is also an inverse Wishart distribution. One can sample $\bSigma$ from that distribution directly. 

\section{Proofs}
\label{sec-proof}
The lower bound for the denominator in the expression for the posterior probability obtained in the following result relies on sufficient prior concentration near the truth and is instrumental in establishing the posterior contraction rate. 
Let $f$ stands for the joint density of $(Y_1,\ldots,Y_n)$ under a generic value of the parameter $(\bbeta,\bSigma)$ and $f_0$ stand for that under the true value $(\bbeta_0,\bSigma_0)$. 

\begin{lemma}
	\label{lemma-3.1}
	For some constant $C_1>0$, $\mathcal{B}_0$ and $\mathcal{H}_0$ are defined in \eqref{assump-true-parameters}, we have
		$\sup\{\mathbb{P}_0(E_n^c): \bbeta_0 \in \mathcal{B}_0,
	\bSigma_0 \in \mathcal{H}_0\} \rightarrow 0$,
	where the set 
	 $E_n = \Big\{\int \int 
	\frac{f}{f_0}
	d\Pi(\bbeta) d\Pi(\bSigma) 
	\geq e^{-C_1 n\epsilon_n^2}\Big \}$. 
\end{lemma}

\begin{proof}
	In view of Lemma 8.10 of \citet{ghos17}, it suffices that 
	\begin{eqnarray}
	&-\log \Pi \left\{(\bbeta,\bSigma): K(f_0, f) \leq n\epsilon_n^2, 
	V(f_0, f) \leq n\epsilon_n^2 \right\}
	\lesssim 
	n \epsilon_n^2,
	\label{eqn-pthm3.1-1}
	\end{eqnarray}
	where $K(f_0,f)$ and $V(f_0,f)$ respectively stand for the average Kullback-Leibler divergence and average Kullback-Leibler variation between $f_0$ and $f$ given by 
	\begin{equation*}
	\begin{split}
	 \frac{1}{n} K(f_0, f) 
	& = 
	\frac{1}{2}\bigg(
	\text{Tr}(\bSigma^{-1} \bSigma_0) - d
	- \log \det(\bSigma^{-1} \bSigma_0)+ \frac{1}{n}\sum_{i=1}^n
	\lVert\bSigma^{-1/2} 
	(\bbeta - \bbeta_0)'
	X_i'\rVert^2
	\bigg),\\
	\frac{1}{n} V(f_0, f)
	& = 
	\frac{1}{2}\Big(
	\text{Tr}\left((\bSigma^{-1}\bSigma_0)^2\right)
	- 2\text{Tr}(\bSigma^{-1}\bSigma_0) + d
	\Big) +
	\frac{1}{n}
	\sum_{i=1}^n \lVert\bSigma_0^{1/2} \bSigma^{-1}
	(\bbeta - \bbeta_0)'X_i'\rVert^2.
	\end{split}
	\end{equation*}
	Define a set of covariance matrices $\mathcal{A}_1$
	and $\mathcal{A}_2$ a set of pairs of regression coefficients and covariance matrices by 
	\begin{align*}
		\mathcal{A}_1 & =
			\Big\{ \bSigma: 
			\text{Tr}(\bSigma^{-1} \bSigma_0) - d
			- \log \det(\bSigma^{-1} \bSigma_0) \leq \epsilon_n^2,\\
			&\qquad\qquad
			\text{Tr}\left((\bSigma^{-1}\bSigma_0)^2\right)
			- 2\text{Tr}(\bSigma^{-1}\bSigma_0) + d
			\leq \epsilon_n^2
			\Big\},\\
		\mathcal{A}_2 
		& = 
			\Big\{ (\bbeta, \bSigma):
			\sum_{i=1}^n \lVert
			\bSigma^{-1/2} (\bbeta - \bbeta_0)' X_i' \rVert^2
			\leq n \epsilon_n^2, \
			\sum_{i=1}^n \lVert \bSigma_0^{1/2} \bSigma^{-1}
			(\bbeta - \bbeta_0)'X_i' \rVert^2 \leq n \epsilon_n^2/2
			\Big\}.
	\end{align*}
	Then a lower bound for the prior probability in \eqref{eqn-pthm3.1-1} can be obtained by lower bounding $\Pi(\mathcal{A}_1)$ and 	$\Pi(\mathcal{A}_2|\mathcal{A}_1)$ separately and multiplying.

	Writing $\bSigma^* = 
	\bSigma_0^{-1/2} \bSigma \bSigma_0^{-1/2}$, $\mathcal{A}_1$ can be written as
	\begin{align*}
	\mathcal{A}_1 &= \left\{ \bSigma:
	\sum_{k=1}^d \left(\eig_k(\bSigma^{*-1}) - 1 - \log \eig_k(\bSigma^{*-1})\right) \leq \epsilon_n^2,\ 
	\sum_{k=1}^d \left( \eig_k(\bSigma^{*-1})-1 \right)^2
	\leq \epsilon_n^2 
	\right\}.
	%	\label{eqn-pthm3.1-4}
	\end{align*}
	By Taylor's expansion
	$\log(x+1) = x-x^2/2 + o(1)$ as $x\rightarrow0$ and since $\epsilon_n\rightarrow 0$, 
	it follows that the second condition in $\mathcal{A}_1$ implies the first, and hence 
	$\mathcal{A}_1 =\big\{ \bSigma: 
	\sum_{k=1}^d ( \eig_k(\bSigma^{*-1})-1 )^2
	\leq \epsilon_n^2\big \}$ for sufficiently large $n$. Since the eigenvalues of $\bSigma_0$ are between $b_1$ and $b_2$ by Assumption~\ref{assump-true-parameters}, Lemma A.1 of \citet{bane15} gives that $\sum_{k=1}^d ( \eig_k(\bSigma^{*-1})-1 )^2\le b_2^2 \lVert\bSigma^{-1}-\bSigma_0^{-1}\rVert_F^2$, and hence  $\mathcal{A}_1 \supset\big\{ \bSigma:  \lVert \bSigma^{-1}-\bSigma_0^{-1} \rVert_F\le \epsilon_n/b_2\big \}$ for sufficiently large $n$.
	Writing in terms of the  eigendecomposition $\bSigma=\bP\bD\bP'$, the triangle inequality, the norm-inequality $\|\bm{A} \bm{B} \|_F \le \min \{\|\bm{A}\| \| \bm{B}\|_F, \|\bm{A}\|_F \| \bm{B}\|\}$ and the facts that $\|\bm{P}\|=1=\|\bm{P}_0\|$ and $\|\bm{D}_0^{-1}\|$ is bounded, we have that
	\begin{align*}
	\lVert \bSigma^{-1}-\bSigma_0^{-1} \rVert_F
	& \le \lVert\bP_0 \rVert\lVert\bP \rVert \lVert  \bD^{-1}- \bD_0^{-1} \rVert_F  + (\lVert \bP_0\rVert\lVert \bD_0^{-1}\rVert + \lVert \bP\rVert\lVert \bD^{-1}\rVert ) \lVert  \bP-\bP_0  \rVert_F\\
	&\lesssim  \lVert  \bD^{-1}- \bD_0^{-1} \rVert_F  + \lVert  \bP-\bP_0  \rVert_F + \lVert  \bD^{-1}- \bD_0^{-1} \rVert_F  \lVert  \bP-\bP_0  \rVert_F,
	\end{align*}
	since $\lVert \bD^{-1}\rVert\le  \lVert \bD_0^{-1}\rVert + \|\bD^{-1}-  \bD_0^{-1}\|$, and the spectral norm is always bounded by the Frobenius norm. Therefore, we have that
	\begin{align*}
	\mathcal{A}_1 \supset  \left\{ \bSigma:  \lVert\bD^{-1}- \bD_0^{-1} \rVert_F \le \epsilon_n/c_1  , \, \lVert  \bP-\bP_0  \rVert_F \le \epsilon_n/c_1
	\right\},
	\end{align*}	
	for some $c_1>0$. Using the prior independence of the eigenvalue distribution and positive lower bound for the prior density at all concerned true value $\bSigma_0$, it is easy to see that $\log \Pi\left\{ \bSigma:  \lVert\bD^{-1}- \bD_0^{-1} \rVert_F\le \epsilon_n/c_1 \right\} \gtrsim -d\log (1/\epsilon_n) \gtrsim -d\log n$. To lower bound $\Pi (\bP: \lVert  \bP-\bP_0  \rVert_F \le \epsilon_n/c_1)$, note that 
	$\Pi$ is the Haar measure on a compact Lie group of dimension $d(d-1)/2$. This means that all translates of  $\{\bP: \lVert  \bP-\bP_0  \rVert_F \le \epsilon_n/c_1\}$ have the same probability, and $N$ many such translates can cover the entire set of $d\times d$ orthogonal matrices, where $N$ stands for the $\epsilon_n/c_1$-covering number of the set of $d\times d$ orthogonal matrices in terms of the Frobenius distance. A crude upper bound for $N$ is easily obtained by embedding the set of $d\times d$ orthogonal matrices in $[-1,1]^{d^2}$, giving the estimate $N\le (2c_1/\epsilon_n)^{d^2}$. This leads to the estimate $\log \Pi\left\{ \bSigma:  \lVert\bP- \bP_0 \rVert_F\le \epsilon_n/c_1 \right\} \gtrsim -d^2\log (2c_1/\epsilon_n) \gtrsim -d^2\log n$. Thus $\log \Pi(\mathcal{A}_1) \gtrsim -d^2\log n$ using the prior independence of $\bD$ and $\bP$. 
	
	To derive a lower bound for 
	$\Pi(\mathcal{A}_2|\mathcal{A}_1)$,
	we first note that $\|\bSigma^{-1}-\bSigma_0^{-1}\|_F\lesssim \epsilon_n$ implies that $\|\bSigma^{-1}\|$ and $\|{\bSigma^*}^{-1}\|$ are bounded by a fixed constant, and hence 
	$n^{-1}\sum_{i=1}^n X_i
	(\bbeta - \bbeta_0)\bSigma^{-1}
	(\bbeta - \bbeta_0)'X_i'
	$ and  
	$n^{-1}\lVert \bSigma^{*-1}\rVert\lVert\bSigma_0^{-1}\rVert\sum_{i=1}^n \|X_i
	(\bbeta - \bbeta_0)\|^2$ are both bounded by a constant multiple of $n^{-1} \sum_{i=1}^n \|X_i
	(\bbeta - \bbeta_0)\|^2= 	n^{-1}\lVert \bX(\bbeta-\bbeta_0) \rVert_F^2 $. 
	Now by the inequality
	\begin{equation}
	\label{normbound}
	\lVert \bX(\bbeta-\bbeta_0) \rVert_F \le \lVert \bX\rVert_\circ\sum_{j=1}^G\lVert \bbeta_j-\bbeta_{0j}\rVert_F\le \lVert \bX\rVert_\circ \sum_{k=1}^d \lVert\beta_k-\beta_{0,k}\rVert_{2,1},
	\end{equation}
	to bound 
	$\Pi(\mathcal{A}_2|\mathcal{A}_1)$ from below, it suffices to bound 
	$\Pi\left( \sum_{k=1}^d\lVert\beta_k-\beta_{0,k}\rVert_{2,1}
	\leq c r_n
	\right)$, where 
	$r_n = 
	\sqrt{{n\epsilon_n^2}}/{ \|\bX\|_\circ}$ and $c$ is a positive constant.
	By (\ref{eqn-prior-beta}), this can be further bounded below by
	\begin{align}
	\pi(s_0)\frac{1}{{Gd \choose s_0}}    
	\int_{\sum_{k=1}^d\|\beta_{S_{0,k}} - \beta_{0,S_{0,k}}\|_{2,1} 
		\leq 
		{c r_n }}
	\prod_{k=1}^d	g_{s_{0,k}}(\beta_{S_{0,k}}) d\beta_{S_{0,1}}\dots d\beta_{S_{0,d}}.
	\label{eqn-pthm3.1-7}
	\end{align}
	By changing the variable 
	$\beta_{S_{0,k}} - \beta_{0, S_{0,k}}$ to $\check{\beta}_{S_{0,k}}$
	and using the fact that $\|x\| \leq \|x\|_1$ 
	for any vector $x$,
	the integral in (\ref{eqn-pthm3.1-7}) is bounded below by
	\begin{align*}
	& e^{-\sum_{k=1}^d\lambda_k\|\beta_{0,k}\|_{2,1}}
	\int_{\sum_{k=1}^d\|\check{\beta}_{S_{0,k}}\|_1 
		\leq 
		{c r_n }
	}
	\prod_{k=1}^d	g_{s_{0,k}}(\check\beta_{S_{0,k}})
	d\check{\beta}_{S_{0,1}}\dots d\check{\beta}_{S_{0,d}}\\
	& \quad\ge e^{-\sum_{k=1}^d\lambda_k\|\beta_{0,k}\|_{2,1}}
	\prod_{k=1}^d
	\prod_{j\in S_{0,k}}
	\left(\frac{2\lambda_k}{ a_j \overline \lambda}\right)^{p_j}\\
	& \qquad	\times 
	\int_{\sum_{k=1}^d\|\check{\beta}_{S_{0,k}}\|_1 
		\leq 
		{c r_n}
	}
	\left(\frac{\overline \lambda}{2}\right)^{p_{S_0}}
	e^{
		-\overline\lambda\sum_{k=1}^d\|\check{\beta}_{S_{0,k}}\|_1 
	}
	d\check{\beta}_{S_{0,1}}\dots d\check{\beta}_{S_{0,d}}.
	\end{align*}
	Using 
	the result that
	the integrand equals to the probability of the first 
	$p_{S_0}$ events of a Poisson process happen before 
	time ${cr_n }$ (similar to the argument used to derive (6.2) in \citealp{cast15}), the last display is further bounded below by
	\begin{align*}
	&e^{-\sum_{k=1}^d\lambda_k \|\beta_{0,k}\|_{2,1}}\Bigg\{	\prod_{k=1}^d
	\prod_{j\in S_{0,k}}
	\left(\frac{2\lambda_k}{ a_j \overline \lambda}\right)^{p_j}\Bigg\}
	e^{-\overline\lambda cr_n}
	\frac{1}{p_{S_0}!}
	\left(
	{\overline\lambda c r_n}\right)^{p_{S_0}}\\
	&\quad\ge e^{-\sum_{k=1}^d\lambda_k\|\beta_{0,k}\|_{2,1}-\overline\lambda c r_n}
	\Bigg\{	\prod_{k=1}^d
	\prod_{j\in S_{0,k}}
	\left(\frac{2}{ a_j }\right)^{p_j}\Bigg\}
	\frac{1}{p_{S_0}!}
	\left(
	{\underline\lambda c r_n}\right)^{p_{S_0}}.
	\end{align*}
	Hence, by Assumption \ref{assump-prior-dimension}, (\ref{eqn-pthm3.1-7}) is bounded below by
	\begin{align*}
	\frac{\pi(0) A_1^{s_0}}{(G\vee n^{p_{\max}})^{A_3s_0} (Gd)^{s_0}}
	e^{- \sum_{k=1}^d \lambda_k\|\beta_{0,k}\|_{2,1}-\overline\lambda c r_n}
	\frac{\left(
		{\underline\lambda c r_n}\right)^{p_{S_0}}}{p_{S_0}!}
	\prod_{k=1}^d
	\prod_{j\in S_{0,k}}
	\left(\frac{2}{ a_j }\right)^{p_j}, 
	\end{align*}
	implying that 
	$\log \Pi(K(f_0, f) \leq n\epsilon_n^2, V(f_0, f) \leq n\epsilon_n^2)$
	is bounded below by
	\begin{align}
	& -d^2\log n +\log \pi(0) + s_0 \log A_1 
	-c_{14} s_0 (\log G+p_{\max}\log n+\log d)-\sum_{k=1}^d  \lambda_k\|\beta_{0,k}\|_{2,1}
	\nonumber\\
	& \quad	
	- \overline\lambda c r_n + p_{S_0}\log(\underline\lambda c r_n)-\log(p_{S_0}!) - \sum_{k=1}^d \sum_{j\in S_{0,k}} p_j \log ( a_j/2),
	\label{eqn-pthm3.1-9}
	\end{align}
	for some constant $c_{14}>0$.
	As $\pi(0)$ is bounded away from zero,
	and Assumption \ref{assump-lambda} gives 
	$
	\overline\lambda c r_n
	- p_{S_0}
	\log \left(\underline\lambda c r_n \right)
	\lesssim \sqrt{n}\epsilon_n\sqrt{\log G} + ({p_{S_0}}/{p_{\max}})\log G  \lesssim n\epsilon_n^2$, the second, sixth, and seventh terms are controlled.
	
	Also, since
	$\sum_{k=1}^d  \|\beta_{0,k}\|_{2,1} \leq \overline{\beta}$
	with 
	the expression of $\overline{\beta}$ is displayed in (\ref{eqn-true-parameters}),
	we have
	$\sum_{k=1}^d \lambda_k  \|\beta_{0,k}\|_{2,1} \le\max_{1 \leq k \leq d} 
	\lambda_k\sum_{k=1}^d   \|\beta_{0,k}\|_{2,1} \leq n\epsilon_n^2$.
	Furthermore,
	since $\log (p_{S_0}!) \leq p_{S_0} \log p_{S_0}$
	and $a_j = O(p_j^{1/2})$, we obtain that
	$
	\log(p_{S_0}!) + \sum_{k=1}^d \sum_{j\in S_{0,k}} p_j \log ( a_j/2) \lesssim s_0 p_{\max} \log n\le n\epsilon_n^2.
	$
	Thus  \eqref{eqn-pthm3.1-9} is bounded below by a constant multiple of $-n\epsilon_n^2$. 
\end{proof}

\begin{proof}[Proof of Lemma \ref{lemma-dim}]
	Let
	$\mathcal{B}_n 
	= 
	\{\bbeta: s_{\bbeta}< r \}$.
	We show that 
	$\mathbb{E}_0\Pi(\bbeta\in\mathcal{B}_n^c | Y_1, \dots, Y_n)
	\rightarrow 0$ as $n \rightarrow \infty$  for $r \geq s_0$. 
	By Lemma \ref{lemma-3.1}, the denominator of (\ref{eqn-post-mixture}) in the expression for $\Pi(\bbeta\in\mathcal{B}_n^c | Y_1, \dots, Y_n)$ with $\mathcal{B}_n$ as above,  
	is bounded below by   $e^{-C_1 n\epsilon_n^2}$  with a large probability. 
	To derive an upper bound for the corresponding numerator, 
	note that its expected value is 
	\begin{align*}
	%\label{eqn-plemma4.1-2}
	\mathbb{E}_0\Big(
	\int\int_{\mathcal{B}_n^c} 
	({f}/{f_{0}})
	d\Pi(\bbeta)d\Pi(\bSigma)
	\Big)
	\leq 
	\int_{\mathcal{B}_n^c}
	d\Pi(\bbeta) 
	=
	\Pi(s_{\bbeta} \geq r) =	\sum_{s = r}^\infty 
	\pi(s),
	\end{align*}
	and by Assumption \ref{assump-prior-dimension}
	and $A_2/(G\vee n^{p_{\max}})^{A_4} \leq 1/2$ as $n \rightarrow \infty$, the bound simplifies to 
	$$\pi(s_0)\Big(\frac{A_2}{(G\vee n^{p_{\max}})^{A_4}} \Big)^{r-s_0}
	\sum_{j = 0}^\infty 
	\Big(\frac{A_2}{(G\vee n^{p_{\max}})^{A_4}} \Big)^j \\
	\leq 
	2 \Big(\frac{A_2}{(G\vee n^{p_{\max}})^{A_4}} \Big)^{r-s_0}.
	$$
	Therefore, because $\mathbb{E}_0\Pi(\mathcal{B}_n^c|Y_1, \dots, Y_n) \leq 
	\mathbb{E}_0\Pi(\mathcal{B}_n^c|Y_1, \dots, Y_n)
	\mathbbm{1}_{E_n}
	+ \mathbb{P}_0(E_n^c)$ and $\mathbb{P}_0(E_n^c)\to 0$, choosing $r=M_2\{s_0\vee [d^2\log n / (\log G\vee p_{\max}\log n)]\}$ for some $M_2$ large enough, we obtain that $\mathbb{E}_0\Pi(\mathcal{B}_n^c|Y_1, \dots, Y_n)$ is bounded above by 
	$$
	\exp\Big( C_1n\epsilon_n^2 
	+ \log 2 + 
	(r-s_0) (\log A_2 - A_4 (\log G\vee p_{\max}\log n))
	\Big) + o(1)\to 0.$$
\end{proof}

\begin{proof}[Proof of Theorem \ref{thm-3.1}]
	
	The proof contains two parts.
	In the first part, we obtain the posterior contraction rate with respect to the average negative log-affinity.
	In the second part, we use the results obtained from the first part to derive (\ref{eqn-post-beta}) and (\ref{eqn-post-sigma}).
	
	\bigskip
	{\bf Part I.}
	Note that for every $\epsilon>0$,
	\begin{align*}
	&\mathbb{E}_0\Pi\left( (\bbeta,\bSigma)\in\mathbb{R}^{p\times d}\times\mathcal H : \frac{1}{n}\sum_{i=1}^n \rho( f_i,f_{0,i}) > \epsilon|Y_1, \dots, Y_n\right)\\
	&\quad\le\mathbb{E}_0\Pi\left((\bbeta,\bSigma)\in\mathcal{B}_n\times\mathcal H : \frac{1}{n}\sum_{i=1}^n \rho( f_i,f_{0,i}) > \epsilon|Y_1, \dots, Y_n\right)
	+\mathbb{E}_0\Pi(\mathcal{B}_n^c|Y_1, \dots, Y_n),
	\end{align*}
	where $\mathcal H$ is the space of $d\times d$ positive definite matrices and $\mathcal{B}_n=\{\bbeta: s_{\bbeta}< M_2 s^\star\}$. The second term on the right hand side goes to zero by Lemma \ref{lemma-dim}, and hence it suffices to show that the first term goes to zero  for  $\epsilon^2=M_1 \epsilon_n^2$.

	Define the sieve
	\begin{align*}
	\mathcal{F}_n = 
	\Bigg\{
	(\bbeta,\bSigma)\in \mathcal{B}_n\times\mathcal H:  
	\max_{\substack{1\leq j\leq G \\
			1 \leq k \leq d}}
	\|\beta_{jk}\|
	\le  H_n, \,
	n^{-1}< \eig_1({\bSigma}^{-1}), \,
	\eig_d({\bSigma}^{-1}) 
	\le n
	\Bigg\},
	\end{align*}
	where $	H_n = {p_{\max}n}/	{\underline\lambda}$ for  $\underline{\lambda}$ given in (\ref{eqn-assump-2}).
	Then 
	\begin{equation}
	\begin{split}
	\Pi( (\mathcal{B}_n\times\mathcal H)\setminus \mathcal{F}_n )
	& \leq
	\sum_{S: s \leq M_2 s^\star} \frac{\pi(s)}{\binom{Gd}{s}}
	\sum_{k=1}^d 
	\sum_{j \in S_k}
	\Pi(\|\beta_{jk}\| \ge H_n)
	\\
	& \quad +
	\Pi\left(\eig_1({\bSigma}^{-1}) \leq n^{-1}\right) 
	+ \Pi\left(\eig_d({\bSigma}^{-1}) \geq n\right).
	\label{eqn-pthm3.1-11}
	\end{split}	
	\end{equation}
	It is easy to see that  $\|\beta_{jk}\|$ is gamma distributed with shape parameter $p_j$ and scale parameter $\lambda_k$.  
	Applying the estimate of the tail of a gamma density
	on page 29 of \citet{bouc13} 
	and the inequality
	$1+x-\sqrt{1+2x} \geq (x-1)/2$,
	for any $x > 0$, we have that 
	$$
	\Pi(\|\beta_{jk}\| >  H_n)
	\leq 
	\exp\left(-p_j \left(1 + \frac{\lambda_k H_n}{p_j} -
	\sqrt{1 + 2 \frac{\lambda_k H_n}{p_j}} \right)\right) 
	\leq 
	\exp \left(- \underline\lambda H_n +p_{\max}
	\right),
	%\label{eqn-pthm3.1-10}
	$$
	for $j = 1, \dots, G$, $k = 1, \dots, d$, leading to the estimate 
	$$
	\sum_{s=1}^{M_2 s^\star} \pi(s)
	s	\exp \left(- \underline\lambda H_n +p_{\max}
	\right) \le \exp\Big( \log (M_2 s^\star) -p_{\max}(n-1)\Big).
	$$
	The second and third terms in 
	(\ref{eqn-pthm3.1-11}) are both bounded by $e^{-c_2 n}$  for some $c_2>0$ by the tail property of inverse Gaussian distribution.
	Combining all these estimates, we obtain that 
	for all sufficiently large $n$,
	\begin{align*}
	\Pi( (\mathcal{B}_n\times\mathcal H) \setminus\mathcal{F}_n) 
	\leq \exp\left(-(1+C_1)n \epsilon_n^2\right).
%	\label{eqn-pthm3.1-2}
	\end{align*}
	
	Next, we construct a test  $\varphi_n$ such that 
	\begin{align}
	\label{eqn-pthm3.1-3}
	\mathbb{E}_{f_0}\varphi_n \lesssim e^{-M_1 n\epsilon_n^2/2},
	\quad 
	\sup_{
		\substack{f \in \mathcal{F}_n: 
			\rho(f_0, f) > M_1  n\epsilon_n^2}}
	\mathbb{E}_{f}(1 - \varphi_n) \lesssim  e^{-M_1 n\epsilon_n^2},
	\end{align}
	for some $M_1>C_1+1$, where 
	$f_0 = \prod_{i=1}^n f_{0,i}$, $f_{0,i} = \mathcal{N}(X_i\bbeta_0, \bSigma_0)$ and 
	$f = \prod_{i=1}^n f_{i}$,  $f_{i} = \mathcal{N}(X_i\bbeta, \bSigma)$, $i=1,\ldots,n$, as required for an application of the general theory of posterior contraction. To this end, we first consider testing $H_0:f=f_0$ against a single point $f=f_1$ in the alternative. Consider the most powerful Neyman-Pearson test 
	$\phi_n = \mathbbm{1}\{f_1/f_0 \geq 1 \}$.
	If the average R\'eyni divergence 
	$-n^{-1}\log \int f_0^{1/2} f_1^{1/2}$
	between $f_0$ and $f_1$
	is bigger than $\epsilon^2>0$,
	then 
	\begin{align*}
	\mathbb{E}_{f_0}\phi_n 
	= \mathbb{E}_{f_0}\left(\sqrt{{f_1}/{f_0}} \geq 1\right)
	\leq \int \sqrt{f_0 f_1} \leq e^{-n\epsilon^2}, \\ 	\mathbb{E}_{f_1}(1-\phi_n) 
	= \mathbb{E}_{f_1}\left(\sqrt{{f_0}/{f_1}} \geq 1\right)
	\leq \int\sqrt{f_0 f_1} \leq e^{-n\epsilon^2}.
	%\label{eqn:test1}
	\end{align*}
	The test $\phi_n$ can also have exponentially small probability of type II error at other alternatives, because by the Cauchy-Schwarz inequality, 
	\begin{align}
	\mathbb{E}_f(1 - \phi_n) 
	\leq 
	\left\{\mathbb{E}_{f_1} (1-\phi_n)\right\}^{{1}/{2}}
	\big\{\mathbb{E}_{f_1}
	\left({f}/{f_1}\right)^{2}\big\}^{{1}/{2}}.
	\label{eqn:test2}
	\end{align}
	so that the expression can be controlled properly if the second factor grows at most like $e^{c n\epsilon^2}$ where $c>0$ can be chosen suitably small. 
	Now we show that
	$\mathbb{E}_{f_1}({f}/{f_1})^{2}$ is bounded for every density with parameters such that
	\begin{align}
	\|\bbeta_1 - \bbeta\|_\infty \leq \frac{1}{   s^\star \sqrt{p_{\max}n} \lVert \bX \rVert_\circ}, \quad 
	\|\bSigma_1 - \bSigma\| \leq   \frac{1}{n^2d},\quad \|\bSigma^{-1}\|\le n.
	\label{eqn:smallpiece}
	\end{align}
	To see this, we observe that for
	$\bSigma^{\star}_1 = \bSigma^{-1/2} \bSigma_1 \bSigma^{-1/2}$,
	\begin{equation}
	\begin{split}
	\mathbb{E}_{f_1} ({f}/{f_1})^{2}
	& = \left(\det(\bSigma^\star_1)\right)^{n/2}
	\left(\det(2\bI - {\bSigma^\star_1}^{-1})\right)^{-n/2}\\
	& \quad 
	\times \exp\Big(\sum_{i=1}^n
	X_i (\bbeta - \bbeta_1) \bSigma^{-1/2}
	(2\bSigma^\star_1 - {\bI})^{-1}
	\bSigma^{-1/2} (\bbeta - \bbeta_1)' X_i'
	\Big).
	\label{eqn-pthm3.1-13}
	\end{split}    
	\end{equation}
	Because $\bSigma\in \mathcal{F}_n$, the condition 
	$\|\bSigma_1 - \bSigma\| \leq \delta_n' = 1/(n^2d)$ implies
	that 
	\begin{align*}
	\|\bSigma^\star_1 - \bI\| \leq 
	\|\bSigma^{-1}\|\|\bSigma_1 - \bSigma\|
	\leq n\|\bSigma_1 - \bSigma\| \leq n\delta_n',
	\end{align*}
	and hence
	$1 - n\delta_n' \leq \eig_1(\bSigma^\star_1)\leq \eig_d(\bSigma^\star_1) \leq 1 + n\delta_n'$.
	Therefore, we obtain that
	\begin{align*}
	\left(\frac{\det(\bSigma^\star_1)}{\det(2\bI - {\bSigma^\star_1}^{-1})}\right)^{n/2}
	& = \exp\left(
	\frac{n}{2}\sum_{k=1}^d \log \left(\eig_k(\bSigma^\star_1)\right) -
	\frac{n}{2}\sum_{k=1}^d \log \left(2 - \frac{1}{\eig_k(\bSigma^\star_1)}\right)
	\right)\\
	& \leq 
	\exp\left(
	\frac{dn}{2} \log (1+n\delta_n')
	- \frac{dn}{2} 
	\log \left(1-\frac{n\delta_n'}{1 - n\delta_n'}\right)
	\right).
	\end{align*}
	By the inequalities
	$1 - x^{-1} \leq \log x \leq x - 1$ for $x > 0$, the display is further bounded by
	\begin{align*}
	\exp\left(	\frac{n^2d \delta_n'}{2} + 	\frac{d n }{2}\left(\frac{n\delta_n'}{1-2n\delta_n'}\right) \right) \le \exp\left( n^2d \delta_n'\right) = e.
	\end{align*}
	
	By the inequality \eqref{normbound},
	we bound the exponential term in (\ref{eqn-pthm3.1-13}) by
	\begin{align*}
	&	\|\bSigma^{-1}\| \|(2\bSigma^\star_1 - \bI)^{-1}\|
	\sum_{i=1}^n \|X_i (\bbeta_1 - \bbeta)\|^2_2\\
	&\quad \leq 
	\|\bSigma^{-1}\| \
	\|(2\bSigma^\star_1 - \bI)^{-1}\|
	\lVert \bX\rVert_\circ^2\Big(\sum_{k=1}^d\lVert\beta_{1,k}-\beta_{k}\rVert_{2,1}\Big)^2.
	\end{align*}
	Since 
	$\displaystyle
	\|(2\bSigma^\star_1 - \bI)^{-1}\| \leq  2$, 
	$\|\bSigma^{{-1}}\| \le n$, and $\sum_{k=1}^d\lVert\beta_{1,k}-\beta_{k}\rVert_{2,1}\le s_{\bbeta_1-\bbeta} \sqrt{p_{\max}} \lVert \bbeta_1-\bbeta \rVert_\infty\le 2M_2 s^\star \sqrt{p_{\max}} \lVert \bbeta_1-\bbeta \rVert_\infty$ on ${\cal F}_n$,
	the display is further bounded by
	\begin{align*}	
	 8 M_2^2  n s^{\star 2} p_{\max} \lVert \bX \rVert_\circ^2 \lVert \bbeta_1-\bbeta   \rVert_\infty^2 \le 8 M_2^2.
	\end{align*}
	Hence we conclude that \eqref{eqn:test2} is bounded by a multiple of $e^{-n\epsilon^2}$ for every density with a parameter in the piece. 
	
	The desired test $\varphi_n$ satisfying \eqref{eqn-pthm3.1-3} is obtained as the maximum of all tests $\phi_n$ described above, for each piece required to cover the sieve. To complete the proof of \eqref{eqn-pthm3.1-3}, we need to show that $\log N_\ast \lesssim n\epsilon_n^2$, where $N_\ast$ is
	the number of pieces satisfying \eqref{eqn:smallpiece} needed to cover the sieve $\mathcal F_n$ (see Lemma D.3 of \citet{ghos17}). It is easy to see that $\log N_\ast$ is bounded by
	\begin{align*}
	&\log N\Big(\frac{1} { s^\star \sqrt{p_{\max}n} \lVert \bX \rVert_\circ}, \big\{\bbeta : s_{\bbeta}\le M_2  s^\star,	\max_{\substack{1\leq j\leq G \\1 \leq k \leq d}}
	\|\beta_{jk}\|
	<  H_n\big\}, \lVert \cdot \rVert_\infty \Big)\\
	&\quad+ \log N \Big( \frac{1}{n^2 d} ,\left\{ \bSigma: 	n^{-1} < \eig_1({\bSigma}^{-1}), \
	\eig_d({\bSigma}^{-1}) 
	< n\right\}, \lVert\cdot\rVert\Big).
	\end{align*}
	The first term of the display is bounded by
	\begin{align}	
	&\log N\Big(\frac{1} {  s^\star \sqrt{p_{\max}n} \lVert \bX \rVert_\circ}, \left\{\bbeta :	s_{\bbeta}\le M_2  s^\star, \lVert\bbeta-\bbeta_0\rVert_\infty
	<  H_n\right\}, \lVert \cdot \rVert_\infty \Big) \nonumber\\
	&\quad\le\log \left\{ \binom{Gd}{M_2 s^\star}\Big( 3 \sqrt{p_{\max}n}  s^\star H_n \lVert \bX\rVert_\circ \Big)^{M_2  s^\star_n p_{\max}} \right\} \nonumber \\
	&\quad\lesssim  s^\star \log G +   s^\star p_{\max} (\log n + \log(H_n\lVert \bX\rVert_\circ )
	\label{eqn:entrobeta}
	\end{align}
	while the second term is bounded by
	\begin{align*}	
	\log N \Big( \frac{1}{n^2 d} ,\left\{ \bSigma: 	n^{-1} < \eig_1({\bSigma}^{-1}) \right\}, \lVert\cdot\rVert\Big)
	& \le 	\log N \Big( \frac{1}{n^2 d} ,\left\{ \bSigma: 	\lVert\bSigma\rVert_F < n\sqrt{d}  \right\}, \lVert\cdot\rVert_F\Big)\\
	&\le {d^2} \log  \big( n^3 d^{3/2}   \big),
	\end{align*}
	both of which are bounded by a constant multiple of $n\epsilon_n^2$. 
	
	Choosing $\epsilon=M_1\epsilon_n^2$ for a sufficiently large $M_1>1+C_1$, we thus have \eqref{eqn-pthm3.1-3}.
	We finally obtain that the posterior
	$\Pi\left( \sum_{i=1}^n \rho(f_i, f_{0,i}) > M_1 n\epsilon_n^2|Y_1, \dots, Y_n \right)$ goes to zero in $\mathbb{P}_0$-probability.
	
	\bigskip
	
	{\bf Part II.}
	Observe that
	$n^{-1} \sum_{i=1}^n
	\rho(f_i, f_{0,i})$ is equal to  
	$$
	-\log \left(
	\frac{\left(\det(\bSigma)\right)^{1/4}
		\left(\det(\bSigma_0)\right)^{1/4}}
	{\left(\det\left((\bSigma + \bSigma_0)/{2}\right)\right)^{1/2}}
	\right) + \frac{1}{8n} \sum_{i=1}^n
	X_i (\bbeta - \bbeta_0)
	\left(\frac{\bSigma + \bSigma_0}{2}\right)^{-1}
	(\bbeta - \bbeta_0)' X_i'.
	$$
	Then $\sum_{i=1}^n\rho(f_i, f_{0,i}) \lesssim n\epsilon_n^2$
	implies the relations 
	\begin{eqnarray}
	\label{eqn-pthm3.1-14}
	-\log \left(
	\frac{\left(\det(\bSigma)\right)^{1/4}
		\left(\det(\bSigma_0)\right)^{1/4}}
	{\left(\det\left(({\bSigma + \bSigma_0})/{2}\right)\right)^{1/2}}
	\right)
	\lesssim \epsilon_n^2,
	\\
	\label{eqn-pthm3.1-15}
	\frac{1}{8n} \sum_{i=1}^n
	X_i (\bbeta - \bbeta_0) 
	\left((\bSigma + \bSigma_0)/{2}\right)^{-1}
	(\bbeta - \bbeta_0)' X_i'
	\lesssim \epsilon_n^2.
	\end{eqnarray}
	First, we show that the probability of (\ref{eqn-pthm3.1-14}) goes to 1 implies (\ref{eqn-post-sigma}).
	Let 
	\begin{equation*}
	d^2(\bSigma, \bSigma_0) 
	= h^2\left(\mathcal{N}(\boldsymbol 0, \bSigma),
	\mathcal{N}(\boldsymbol 0, \bSigma_0)\right)
	= 1 - 
	\frac{\left(\det(\bSigma)\right)^{1/4}
		\left(\det(\bSigma_0)\right)^{1/4}}
	{\left(\det\left((\bSigma + \bSigma_0)/{2}\right)\right)^{1/2}}.
	\end{equation*}
	Since the eigenvalues of $\bSigma_0$ lie in $[b_1,b_2]$, by Lemma 2 of \citet{suar17}, we obtain that
	$
	%\label{eqn-pthm3.1-16}
	d^2(\bSigma, \bSigma_0) 
	\gtrsim
	\| \bSigma_0^{-1/2} (\bSigma - \bSigma_0)
	\bSigma_0^{-1/2}\|_F^2,
	$
	if the left hand side is sufficiently small.
	Since 
	\begin{align*}
	-\log \left(
	\frac{\left(\det(\bSigma)\right)^{1/4}
		\left(\det(\bSigma_0)\right)^{1/4}}
	{\left(\det\left((\bSigma + \bSigma_0)/{2}\right)\right)^{1/2}}
	\right)
	=
	-\log(1-d^2(\bSigma, \bSigma_0))
	\ge
	d^2(\bSigma, \bSigma_0),
	\end{align*}
	we obtain that
	$\| \bSigma - \bSigma_0 \|_F^2 \lesssim \epsilon_n^2$. This proves (\ref{eqn-post-sigma}).
	
	Next, we show that the probability (\ref{eqn-pthm3.1-15}) goes to 1 implies (\ref{eqn-post-beta}).
	Given (\ref{eqn-post-sigma})
	and by Assumption \ref{assump-true-parameters}, we obtain that
	\begin{align*}
	\|\bSigma + \bSigma_0\|^2 = 
	\|\bSigma - \bSigma_0 + 2\bSigma_0\|^2
	\leq 
	2\|\bSigma - \bSigma_0\|_F^2 + 8\|\bSigma_0\|^2
	\lesssim \epsilon_n^2 + 1.
	\end{align*}
	Hence using 
	$\eig_1 \left(({\bSigma+\bSigma_0}/{2})^{-1}\right) 
	= \left(\eig_d ({\bSigma+\bSigma_0}/{2})\right)^{-1} 
	=\left\|({\bSigma+\bSigma_0}/{2})\right\|^{-1}\ge (1+\epsilon_n^2)^{-1/2}
	$, (\ref{eqn-pthm3.1-15}) implies that
	\begin{align*}
	\epsilon_n^2 \ge 
	\frac{1}{8n}\sum_{i=1}^n 
	\|X_i (\bbeta - \bbeta_0)\|^2
	\Big\|\frac{\bSigma + \bSigma_0}{2}\Big\|^{-1}
	&\gtrsim
	\frac{1}{n} \sum_{i=1}^n
	\|X_i(\bbeta - \bbeta_0)\|^2 
	/\sqrt{\epsilon_n^2 + 1}.
	\end{align*}
	Combining with (\ref{eqn-post-sigma}),
	we obtain (\ref{eqn-post-beta}).
\end{proof}

%\subsection{Proofs of Theorem \ref{thm-bvm} and \ref{thm-selection}}
\begin{proof}[Proof of Theorem \ref{thm-bvm}]
	Let $\mathcal{H}_n 
	= \{
	\bSigma\in{\cal H}: 
	\|\bSigma - \bSigma_0\|_F^2 \leq 
	M_1\epsilon_n^2
	\}$ and 
	$$
	\Theta_n 
	= \left\{ \bbeta\in\mathbb{R}^{p\times d}:
	s_{\bbeta} \leq M_2s^\star, \bigg(
	\sum_{k=1}^d
	\|\beta_k - \beta_{0,k}\|_{2,1}\bigg)^2
	\le
	\frac{M_4 n\epsilon_n^2 s^\star }
	{
		\lVert \bX\rVert_\circ^2
		\phi_{\ell_{2,1}}^2(s_0+M_2s^\star)
	}
	\right\}, 
	$$
	where ${\cal H}$ is a space of $d\times d$ positive definite matrices.
	The proof contains two parts. 
	In the first part,
	we show that the total variation metric
	between  $\Pi(\bbeta\in \cdot|Y_1, \dots, Y_n)$
	and $\check\Pi_n (\bbeta\in\cdot|Y_1, \dots, Y_n) := \check\Pi_n((\bbeta,\bSigma)\in \cdot\times {\cal H}_n|Y_1, \dots, Y_n) $
	is small,
	where $\check\Pi_n((\bbeta,\bSigma)\in \cdot\times\cdot|Y_1, \dots, Y_n)$ is the renormalized measure of
	$\Pi((\bbeta,\bSigma)\in\cdot\times\cdot|Y_1, \dots, Y_n)$ restricted to the 
	set $\Theta_n\times {\cal H}_n$.
	We also show that the total variation distance 
	between $\Pi^{\infty}(\bbeta\in\cdot|Y_1, \dots, Y_n)$
	and $\check\Pi_n^{\infty}(\bbeta\in\cdot|Y_1, \dots, Y_n)$ 
	is small,
	where $\check\Pi_n^{\infty}(\bbeta\in\cdot|Y_1, \dots, Y_n)$ is the measure $\Pi^{\infty}(\bbeta\in\cdot|Y_1, \dots, Y_n)$ restricted and renormalized 
	to $\Theta_n$.
	In the second part, 
	we show that the total variation distance between $\check\Pi_n(\bbeta\in\cdot|Y_1, \dots, Y_n)$ and  $\check\Pi_n^{\infty}(\bbeta\in\cdot|Y_1, \dots, Y_n)$ is small.

	For any set $A$, let $\Pi_A(\cdot)$ be the renormalized measure of $\Pi(\cdot)$ which is restricted to the set $A$.
	Then $\|\Pi(\cdot) - \Pi_A(\cdot)\| \leq 2\Pi(A^c)$.
	Clearly,
	\begin{align*}
	\mathbb{E}_0\|\Pi(\bbeta\in\cdot|Y_1, \dots, Y_n) - 
	\check\Pi_n(\bbeta\in\cdot|Y_1, \dots, Y_n)\|_{TV}\rightarrow 0,
	\end{align*}
	by (\ref{eqn-post-sigma}) and (\ref{eqn-recover-l2l1}).
	To show that
	\begin{align*}
	%\label{eqn-pthm4.8-0}
	\mathbb{E}_0
	\|\Pi^\infty(\bbeta\in\cdot| Y_1, \dots, Y_n)-
	\check\Pi_n^\infty(\bbeta\in\cdot| Y_1, \dots, Y_n)
	\|_{TV} \rightarrow 0,
	\end{align*}
	we write
	\begin{align}
	\label{eqn-pthm4.8-1}
	\Pi^\infty(\bbeta \in \Theta_n^c|Y_1, \dots, Y_n)
	=
	\frac{
		\int_{\Theta_n^c} 
		\exp\{\ell_n (\bbeta, \bSigma_0 )
		- \ell_n (\bbeta_0, \bSigma_0)
		\} dU(\bbeta)
	}
	{
		\int
		\exp\{\ell_n (\bbeta, \bSigma_0 )
		- \ell_n (\bbeta_0, \bSigma_0)
		\} dU(\bbeta)
	},
	\end{align}
	with $ dU ( \bbeta)$ defined in \eqref{eqn:appdu}.
	By (\ref{eqn-likelihood}),
	$\ell_n(\bbeta, \bSigma_0) - \ell_n(\bbeta_0, \bSigma_0)$
	equals to 
	\begin{align*}
	- \frac{1}{2}\sum_{i=1}^n \|
	\bSigma_0^{-1/2}  \tilde \bX_{i}' \vc(\bbeta - \bbeta_0)'\|^2 
	+ \sum_{i=1}^n
	\left(Y_i - \vc(\bbeta_0) \tilde\bX_{i} \right) \bSigma_0^{-1} 
	\tilde\bX_{i}'  \vc(\bbeta - \bbeta_0)'.
	\end{align*}
	%    Note that except for the elements of the first row
	%  of $\tilde\bX_{i,k}$, the rest are zero.
	By plugging-in the last display into (\ref{eqn-pthm4.8-1}),
	the denominator is bounded below by
	\begin{align*}
	& 
	\frac{\pi(s_0)} {{Gd \choose s_0}}\Bigg(\prod_{k=1}^d \prod_{j \in S_{0,k}} 
	\left(\frac{\lambda_k}{a_j}\right)^{p_j}\Bigg)
	\\
	&\quad\times\int 
	\exp\left(-\frac{1}{2}\sum_{i=1}^n \|
	\bSigma_0^{-1/2}  \tilde \bX_{i,S_0}' \tilde\beta_{S_0}\|^2 +\sum_{i=1}^n 
	\left(Y_i - \vc(\bbeta_0) \tilde\bX_{i} \right) \bSigma_0^{-1} 
	\tilde\bX_{i,S_0}'  \tilde \beta_{S_0}
	\right) d \tilde\beta_{S_0},
	\end{align*}
	where $\tilde \beta_{S_0}=\left((\beta_{1,S_{0,1}}-\beta_{0,1,S_{0,1}})',\dots,(\beta_{d,S_{0,d}}-\beta_{0,d,S_{0,d}})'\right)'$.
	By Jensen's inequality, 
	the display is bounded below by
	\begin{align}
	&    \frac{\pi(s_0)}{{Gd \choose s_0}}
	\Bigg(\prod_{k=1}^d \prod_{j \in S_{0,k}}\left(\frac{\lambda_k}{a_j}\right)^{p_j}\Bigg)
	\int 
	\exp\left(-\frac{1}{2}\sum_{i=1}^n \|
	\bSigma_0^{-1/2}  \tilde \bX_{i,S_0}' \tilde\beta_{S_0}\|^2 
	\right) d \tilde\beta_{S_0},
	\nonumber \\
	& \quad =
	\frac{\pi(s_0)}{{Gd \choose s_0}}
	\Bigg(\prod_{k=1}^d \prod_{j \in S_{0,k}}\left(\frac{\lambda_k}{a_j}\right)^{p_j}\Bigg)
	\sqrt{\frac{(2\pi)^{p_{S_0}}}
		{\det\left(\sum_{i=1}^n \tilde\bX_{i,S_0} \bSigma_0^{-1} \tilde\bX_{i,S_0}'\right)}}.
	\label{eqn-pthm4.8-2}
	\end{align}
	Letting $\bGamma_{S_0} = 
	\sum_{i=1}^n \tilde\bX_{i,S_0} \bSigma_0^{-1} \tilde\bX_{i,S_0}'$,
	we apply Jensen's inequality to obtain that
	\begin{align*}
	\det(\bGamma_{S_0})
	\leq 
	\left({\Tr(\bGamma_{S_0})}/{p_{S_0}} 
	\right)^{p_{S_0}}
	& \leq 
	\Big(\max_{l} (\bGamma_{S_0})_{l,l}\Big)^{
		p_{S_0}},
	\end{align*}  
	where $(\bGamma_{S_0})_{l,l}$ is the $l$th diagonal element of $\bGamma_{S_0}$.
	Note that
	\begin{align*}
	\max_{l} (\bGamma_{S_0})_{l,l}\le \frac{1}{b_1} \max_{1\le j\le G}\left\lVert\sum_{i=1}^n \tilde\bX_{i,j}\tilde\bX_{i,j}'\right\rVert = \frac{1}{b_1}  \max_{1\le j\le G} \lVert  \bX_j \rVert^2 =\frac{\lVert \bX\rVert_\circ^2}{b_1},
	\end{align*}  
	where $\tilde\bX_{i,j}=\bI_d\otimes X_{ij}'$,
	and hence \eqref{eqn-pthm4.8-2} is further bounded below by
	\begin{align}
	&\frac{\pi(s_0)}{{Gd \choose s_0}}
	\Bigg(\prod_{k=1}^d \prod_{j \in S_{0,k}}\left(\frac{\lambda_k}{a_j}\right)^{p_j}\Bigg)
	\left(\frac{2b_1\pi}{\lVert \bX \rVert_\circ^2}\right)^{ p_{S_0}/2}\nonumber\\
	&\quad\ge \frac{\pi(s_0)}{(Gd)^{s_0}\prod_{k=1}^d \prod_{j \in S_{0,k}} a_j^{p_j}}  \left(\frac{\sqrt{2b_1\pi}}{B_1(G^{1/p_{\max}}\vee n)^{B_2}}\right)^{p_{S_0}}\nonumber\\
	&\quad\ge \frac{\pi(s_0)}{(Gd)^{s_0}  a_j^{s_0p_{\max}}}  \left(\frac{\sqrt{2b_1\pi}}{B_1(G^{1/p_{\max}}\vee n)^{B_2}}\right)^{ s_0 p_{\max}}.
	\label{eqn:denlb}
	\end{align}  
	We thus obtain a lower bound for the denominator.
	
	The numerator of (\ref{eqn-pthm4.8-1}) can be written as
	\begin{align}
	\begin{split}
	\int_{\Theta_n^c} 
	\Bigg\{&\exp\left(   
	-\frac{1}{2}\sum_{i=1}^n 
	\|\vc(\bbeta - \bbeta_{0})
	\tilde \bX_{i} \bSigma_0^{-1/2}\|^2 \right)\\
	& 
	\times \exp\left(
	\sum_{i=1}^n \left(Y_i - \vc(\bbeta_{0}) \tilde\bX_{i} \right)
	\bSigma_0^{-1}\tilde \bX_{i}'
	\vc(\bbeta - \bbeta_{0})'
	\right)\Bigg\}
	dU(\bbeta).
	\end{split}
	\label{eqn:numub}
	\end{align}
	Note that
	\begin{align}
	&\sum_{i=1}^n \left(Y_i - \vc(\bbeta_{0}) \tilde\bX_{i} \right)
	\bSigma_0^{-1} \tilde\bX_{i}'
	\vc(\bbeta - \bbeta_{0})' \nonumber\\
	&\quad= \sum_{j=1}^G\sum_{i=1}^n \left(Y_i - \vc(\bbeta_{0}) \tilde\bX_{i} \right)
	\bSigma_0^{-1}  \tilde\bX_{i,j}'
	\vc(\bbeta_j - \bbeta_{0,j})' \nonumber\\
	&\quad\le \sum_{j=1}^G\left\lVert\sum_{i=1}^n \left(Y_i - \vc(\bbeta_{0}) \tilde\bX_{i} \right)
	\bSigma_0^{-1}  \tilde\bX_{i,j}'\right\rVert
	\lVert\bbeta_j - \bbeta_{0,j}\rVert_F.
	\label{eqn:upper111}
	\end{align}
	Using the tail inequality for quadratic forms of Gaussian random variables (Proposition~1 of \citet{hsu2012tail}), we obtain for every $t>0$,
	\begin{align*}
	%\label{eqn-pthm4.8-2.5}
	\mathbb{P}
	\Bigg(&
	\max_{1\le j\le G}\Bigg\lVert\sum_{i=1}^n \left(Y_i - \vc(\bbeta_{0}) \tilde\bX_{i} \right)
	\bSigma_0^{-1}  \tilde\bX_{i,j}'\
	\Bigg\rVert^2\\
	&\qquad
	\geq \Tr(\bDel'\bDel) + 2\sqrt{\Tr((\bDel'\bDel)^2)t} +2 \lVert\bDel\rVert^2 t
	\Bigg) \leq  G e^{-t} ,
	\end{align*}
	where     $\bDel=(\tilde\bX_{1,j} \bSigma_0^{-1},\dots,\tilde\bX_{n,j} \bSigma_0^{-1}) \in\mathbb{R}^{p_j \times dn}$. 
	Since  $\Tr(\bDel'\bDel)\le p_j \lVert\bDel\rVert^2$ and $\lVert\bDel\rVert\lesssim \lVert(\tilde\bX_{1,j},\dots,\tilde\bX_{n,j})\rVert=\lVert\bX_j\rVert\le \lVert\bX\rVert_\circ$, choosing $t=2(\log G\vee p_{\max}\log n)$, we obtain
	\begin{align*}
	\mathbb{P}
	\Bigg(&
	\max_{1\le j\le G}\Bigg\lVert\sum_{i=1}^n \left(Y_i - \vc(\bbeta_{0}) \tilde\bX_{i} \right)
	\bSigma_0^{-1}  \tilde\bX_{i,j}'\
	\Bigg\rVert \ge c_4 \lVert \bX\rVert_\circ \sqrt{\log G \vee p_{\max}\log n}
	\Bigg) \leq  \frac{1}{G},
	\end{align*}
	for some $c_4>0$.
	Let  $D_n=c_4 \lVert \bX\rVert_\circ \sqrt{\log G \vee p_{\max}\log n}$.
	Then, with probability tending to one, \eqref{eqn:upper111} is further bounded by
	\begin{align*}
	D_n \sum_{k=1}^d \lVert\beta_k-\beta_{0,k}\rVert_{2,1}& \le\frac{2 D_n\lVert \bX(\bbeta-\bbeta_0)\rVert_F|S_{\bbeta-\bbeta_0}|^{1/2}}{\lVert \bX\rVert_\circ \phi_{\ell_{2,1}}(|S_{\bbeta-\bbeta_0}|)}-  D_n  \sum_{k=1}^d \lVert\beta_k-\beta_{0,k}\rVert_{2,1}\\
	& = \frac{2D_n\sqrt{|S_{\bbeta-\bbeta_0}|\sum_{i=1}^n\lVert \vc(\bbeta-\bbeta_0) \tilde\bX_i\rVert^2}}{\lVert \bX\rVert_\circ \phi_{\ell_{2,1}}(|S_{\bbeta-\bbeta_0}|)}-  D_n  \sum_{k=1}^d \lVert\beta_k-\beta_{0,k}\rVert_{2,1}.
	\end{align*} 
	The display is further bounded by
	\begin{align*}
	&\frac{2b_2D_n\sqrt{|S_{\bbeta-\bbeta_0}|\sum_{i=1}^n\lVert \vc(\bbeta-\bbeta_0) \tilde\bX_i \bSigma_0^{-1}\rVert^2}}{\lVert \bX\rVert_\circ \phi_{\ell_{2,1}}(|S_{\bbeta-\bbeta_0}|)}-  D_n  \sum_{k=1}^d \lVert\beta_k-\beta_{0,k}\rVert_{2,1}\\
	&\quad\le \frac{1}{2}\sum_{i=1}^n\lVert \vc(\bbeta-\bbeta_0) \tilde\bX_i \bSigma_0^{-1}\rVert^2+\frac{2 b_2^2 D_n^2 |S_{\bbeta-\bbeta_0}|}{\lVert \bX\rVert_\circ^2 \phi_{\ell_{2,1}}^2(|S_{\bbeta-\bbeta_0}|)}-  D_n  \sum_{k=1}^d \lVert\beta_k-\beta_{0,k}\rVert_{2,1},
	\end{align*}
	by the Cauchy-Schwarz inequality. Therefore, with probability tending to one, \eqref{eqn:numub} is bounded by
	\begin{align*}
	&\int_{\Theta_n^c} \exp\left(   \frac{2 b_2^2 D_n^2 |S_{\bbeta-\bbeta_0}|}{\lVert \bX\rVert_\circ^2 \phi_{\ell_{2,1}}^2(|S_{\bbeta-\bbeta_0}|)}-  D_n  \sum_{k=1}^d \lVert\beta_k-\beta_{0,k}\rVert_{2,1}
	\right)dU(\bbeta)\\
	&\quad\le \exp\left(   \frac{2 b_2^2 D_n^2 (s_0+M_2 s^\star)}{\lVert \bX\rVert_\circ^2 \phi_{\ell_{2,1}}^2(s_0+M_2 s^\star)} - \frac{D_n\sqrt{M_4 n\epsilon_n^2 s^\star}}{2\lVert\bX\rVert_\circ \phi_{\ell_{2,1}}(s_0+M_2s^\star)}
	\right) \\
	&\qquad\times  \sum_{S : s\le M_2 s^\star}
	\frac{\pi(s)}{{Gd \choose s}}
	\int_{\Theta_n^c}
	\prod_{k=1}^d \Bigg(\prod_{j\in S_k}\left(\frac{\lambda_k}{a_j}\right)^{p_j}\Bigg)
	\exp\left(-\frac{D_n}{2}\lVert\beta_k-\beta_{0,k}\rVert_{2,1}\right)d \beta_{S_k}\otimes\delta_{S_k^c}.
	\end{align*}
	Since $c_4\lambda_k/B_3\le D_n$ for every $k\le d$,
	the last summation is bounded by
	\begin{align*}
	\sum_{S : s\le M_2 s^\star}
	\frac{\pi(s)}{{Gd \choose s}} \left(\frac{2B_3}{c_4}\right)^{p_S}\le 1,
	\end{align*}
	where the inequality holds by making $c_4$ large enough. 
	Now, plug in $D_n$ and combine the display with \eqref{eqn:denlb} to obtain an upper bound of the expectation of \eqref{eqn-pthm4.8-1}. Since $a_j=O(p_j^{1/2})$ and $\pi(s_0)\gtrsim A_1^{s_0}/(G^{A_3}\vee n^{A_5p_{\max}})^{s_0}$, the upper bound goes to zero as long as $M_4$ is chosen sufficiently large.

	For a measurable subset $\mathcal{B}$ of $\mathbb{R}^{p\times d}$, we can write
	\begin{align*}
	&\check\Pi_n(\mathcal{B}|Y_1, \dots, Y_n)\\
	&~\propto \int_{(\mathcal{B}\cap \Theta_n)} \int_{{\cal H}_n} 
	\frac{\exp(\ell_n(\bbeta,\bSigma))}{\exp(\ell_n(\bbeta_0,\bSigma))} \exp\left(-\sum_{k=1}^d\lambda_{k} \|\beta_{k}\|_{2,1}\right)
	\exp( \ell_n(\bbeta_0,\bSigma))d\Pi(\bSigma)  dU(\bbeta),
	\end{align*}
	and
	\begin{align*}
	& \check\Pi_n^{\infty}(\mathcal{B}|Y_1, \dots, Y_n)\\
	&~\propto \int_{(\mathcal{B}\cap \Theta_n)} \int_{{\cal H}_n} 
	\frac{\exp(\ell_n(\bbeta,\bSigma_0))}{\exp(\ell_n(\bbeta_0,\bSigma_0))} \exp\left(-\sum_{k=1}^d\lambda_{k} \|\beta_{0,k}\|_{2,1}\right)
	\exp( \ell_n(\bbeta_0,\bSigma))d\Pi(\bSigma)  dU(\bbeta).
	\end{align*}
	Note that for sequences of measures ($\mu_S$) and ($\nu_S$), 
	\begin{align*}
	\left\|\frac{\sum_S \mu_S}{\|\sum_S\mu_S\|_{TV}} 
	- \frac{\sum_S \nu_S}{\|\sum_S\nu_S\|_{TV}} \right\|_{TV}
	\leq 2\sup_S \left\|1 - \frac{d\nu_S}{d\mu_S}\right\|_{\infty},
	\end{align*}  
	(see e.g. page 2011 of \citet{cast15}). Hence it suffices to show that
	\begin{align*}
	\mathbb{E}_0\sup_{\bbeta\in\Theta_n}\left\lvert 1-    \frac{\int_{{\cal H}_n} 
		\frac{\exp(\ell_n(\bbeta,\bSigma))}{\exp(\ell_n(\bbeta_0,\bSigma))}\exp\left(-\sum_{k=1}^d\lambda_{k} \|\beta_{k}\|_{2,1}\right)\exp(\ell_n(\bbeta_0,\bSigma))
		d\Pi(\bSigma) }{\int_{{\cal H}_n} 
		\frac{\exp(\ell_n(\bbeta,\bSigma_0))}{\exp(\ell_n(\bbeta_0,\bSigma_0))}\exp\left(-\sum_{k=1}^d\lambda_{k} \|\beta_{0,k}\|_{2,1}\right)
		\exp(\ell_n(\bbeta_0,\bSigma))d\Pi(\bSigma) }\right\rvert \rightarrow 0.
	\end{align*}  
	Using the property that
	$ |1-{\int f  }/{\int g  }|\le (1-\inf ({f}/{g}))\vee (\sup({f}/{g})-1)\le \sup |1-{f}/{g}|$, 
	the expression in the last display is bounded by
	\begin{align*}
	& \mathbb{E}_0\sup_{\bbeta\in\Theta_n}\sup_{\bSigma\in{\cal H}_n}\left\lvert 1- 	\exp\left(\tilde\ell_n (\bbeta,\bSigma)-\sum_{k=1}^d\lambda_{k} (\|\beta_{k}\|_{2,1}-\|\beta_{0,k}\|_{2,1})\right)\right\rvert \\
	&\quad\le \mathbb{E}_0\sup_{\bbeta\in\Theta_n}\sup_{\bSigma\in{\cal H}_n}\Bigg\{\left(|\tilde\ell_n (\bbeta,\bSigma)|+\max_{1\le k\le d}\lambda_k \sum_{k=1}^d\|\beta_{k}-\beta_{0,k}\|_{2,1}\right)\\
	&\quad\qquad\qquad\qquad\qquad\times\exp\left(|\tilde\ell_n (\bbeta,\bSigma)|+\max_{1\le k\le d}\lambda_k \sum_{k=1}^d\|\beta_{k}-\beta_{0,k}\|_{2,1}\right)\Bigg\},
	\end{align*}  
	where $\tilde\ell_n (\bbeta,\bSigma)=\ell_n(\bbeta,\bSigma)+\ell_n(\bbeta_0,\bSigma_0)-\ell_n(\bbeta,\bSigma_0)-\ell_n(\bbeta_0,\bSigma)$.
	First, it is easy to see that $\sup\{\lambda_k \sum_{k=1}^d\|\beta_{k}-\beta_{0,k}\|_{2,1}: \bbeta\in\Theta_n,\,  1\le k\le d\}
	\rightarrow 0$ due to the small $\lambda$ regime. To complete the proof, we shall show that 
	\begin{align}
	\mathbb{E}_0\sup_{\bbeta\in\Theta_n}\sup_{\bSigma\in{\cal H}_n} |\tilde\ell_n (\bbeta,\bSigma)|\rightarrow 0.
	\label{eqn:diffll}
	\end{align}  
	It can be easily verified that
	\begin{align*}
	|\tilde\ell_n (\bbeta,\bSigma)|
	&\le\frac{1}{2}\left|\sum_{i=1}^n 
	\vc(\bbeta-\bbeta_0) \tilde\bX_{i}  (\bSigma^{-1}-\bSigma_0^{-1})
	\tilde \bX_{i}' \vc(\bbeta-\bbeta_0)'\right|\\
	&\quad+\left|\sum_{i=1}^n 
	\vc(\bbeta-\bbeta_0) \tilde\bX_{i}  (\bSigma^{-1}-\bSigma_0^{-1})\big(Y_i - \vc(\bbeta_0) \tilde\bX_{i} \big)'\right|.
	\end{align*}  
	First note that
	\begin{align*}
	&\sup_{\bbeta\in\Theta_n}\sup_{\bSigma\in{\cal H}_n}\left|\sum_{i=1}^n 
	\vc(\bbeta-\bbeta_0) \tilde\bX_{i}  (\bSigma^{-1}-\bSigma_0^{-1})
	\tilde\bX_{i}' \vc(\bbeta-\bbeta_0)'\right|\\
	&\quad\le \sup_{\bSigma\in{\cal H}_n}\lVert \bSigma^{-1}-\bSigma_0^{-1} \rVert \sup_{\bbeta\in\Theta_n} \lVert \bX (\bbeta-\bbeta_0)\rVert_F^2\\
	&\quad\lesssim \lVert\bX\rVert_\circ^2  \sup_{\bSigma\in{\cal H}_n}\lVert \bSigma-\bSigma_0 \rVert \sup_{\bbeta\in\Theta_n} \left( \sum_{k=1}^d 
	\lVert \beta_k-\beta_{0,k} \rVert_{2,1} \right)^2,
	\end{align*}  
	where the last inequality holds by \eqref{normbound} and Assumption~\ref{assump-true-parameters} since $\sup\{\lVert \bSigma-\bSigma_0 \rVert: \bSigma\in{\cal H}_n\}$ is small. The rightmost side of the display is bounded by $s^\star n\epsilon_n^3$ which goes to zero by the assumption.
	Similar to \eqref{eqn:upper111}, we also obtain that
	\begin{align}
	&\mathbb{E}_0\sup_{\bbeta\in\Theta_n}\sup_{\bSigma\in{\cal H}_n}\left| 
	\vc(\bbeta-\bbeta_0) \sum_{i=1}^n\tilde\bX_{i}  (\bSigma^{-1}-\bSigma_0^{-1})\big(Y_i - \vc(\bbeta_0) \tilde\bX_{i} \big)'\right| \nonumber\\
	&\quad \le \mathbb{E}_0\sup_{\bbeta\in\Theta_n}\sup_{\bSigma\in{\cal H}_n} \sum_{j=1}^G
	\left\lVert\vc(\bbeta_j-\bbeta_{0,j})\right\rVert_F\left\lVert W_{\bSigma,j}\right\rVert \nonumber\\
	&\quad\le \mathbb{E}_0 \max_{1\le j\le G}\sup_{\bSigma\in{\cal H}_n} \left\lVert W_{\bSigma,j}\right\rVert  \sup_{\bbeta\in\Theta_n} \sum_{k=1}^d\lVert \beta_k-\beta_{0,k} \rVert_{2,1} ,
	\label{eqn:upper22}
	\end{align}  
	where $W_{\bSigma,j}=\sum_{i=1}^n\tilde\bX_{i,j}  (\bSigma^{-1}-\bSigma_0^{-1})\big(Y_i - \vc(\bbeta_0) \tilde\bX_{i} \big)'$. By Lemma 2.2.2 of \citet{vaart1996weak} applied with $\psi_2(x)=\exp(x^2)-1$, we have
	\begin{align*}
	\mathbb{E}_0 \max_{1\le j\le G} \sup_{\bSigma\in{\cal H}_n} \lVert W_{\bSigma,j} \rVert& \le\sqrt{p_{\max} d }\, \mathbb{E}_0 \max_{1\le j\le G}\max_{1\le \ell\le p_j d}  \sup_{\bSigma\in{\cal H}_n} | W_{\bSigma,j,\ell} | \\
	&\le \sqrt{p_{\max}d }\left\lVert \max_{1\le j\le G}\max_{1\le \ell\le p_j d} \sup_{\bSigma\in{\cal H}_n} | W_{\bSigma,j,\ell} |\right\rVert_{\psi_2}\\
	&\lesssim \sqrt{p_{\max}d \log G} \max_{1\le j\le G} \max_{1\le \ell\le p_j d} \left\lVert  \sup_{\bSigma\in{\cal H}_n} | W_{\bSigma,j,\ell} | \right\rVert_{\psi_2},
	\end{align*}  
	where $\lVert\cdot\rVert_{\psi_2}$ denotes the Orlicz norm and $W_{\bSigma,j,\ell}$ is the $\ell$th element of $W_{\bSigma,j}$. By Lemma 2.2.1 of \citet{vaart1996weak}, we have that for every $\bSigma_1,\bSigma_2\in{\cal H}_n$,
	$$
	\lVert  W_{\bSigma_1,j,\ell} -W_{\bSigma_2,j,\ell} \rVert_{\psi_2}\lesssim \sqrt{{\rm Var}(W_{\bSigma_1,j,\ell} -W_{\bSigma_2,j,\ell})}
	\le\lVert \bSigma_0^{1/2}(\bSigma_1^{-1}-\bSigma_2^{-1})\rVert  \lVert \bX_j\rVert, 
	$$ 
	 which is bounded by $\lVert \bSigma_1 -\bSigma_2  \rVert_F \lVert \bX\rVert_\circ$,
	by the relations $\lVert \bSigma_1 -\bSigma_2  \rVert_F\le  \lVert \bSigma_1 -\bSigma_0  \rVert_F+\lVert \bSigma_2-\bSigma_0  \rVert_F\lesssim \epsilon_n$ and  $\bSigma_1^{-1}-\bSigma_2^{-1}= -\bSigma_1^{-1} (\bSigma_1-\bSigma_2) \bSigma_2^{-1}$, and the fact that the eigenvalues of $\bSigma$, and hence also those of $\bSigma$ and $\bSigma_2$, lie between two fixed positive numbers. We see that $W_{\bSigma,j,\ell}$ is a separable Gaussian process as ${\cal H}_n$ is a separable metric space under the Frobenius norm. Hence, by Corollary 2.2.5 of \citet{vaart1996weak}, for any fixed $\bSigma'\in{\cal H}_n$ and some $c_5>0$,
	\begin{align*}
	\left\lVert  \sup_{\bSigma\in{\cal H}_n} | W_{\bSigma,j,\ell} | \right\rVert_{\psi_2} \lesssim \left\lVert    W_{\bSigma',j,\ell} \right\rVert_{\psi_2}+\int_0^{c_5\lVert \bX\rVert_\circ{\rm diam}_j({\cal H}_n)}\sqrt{\log N\left(\frac{\epsilon}{2c_5\lVert \bX\rVert_\circ},{\cal H}_n ,\lVert\cdot\rVert_F\right)}d\epsilon,
	\end{align*} 
	where ${\rm diam}_j({\cal H}_n)=\sup\{\lVert\bSigma_1-\bSigma_2\rVert_F:\bSigma_1,\bSigma_2\in{\cal H}_n\}$. By Lemma 2.2.1 of \citet{vaart1996weak} again, we have that
	\begin{align*}
	\left\lVert    W_{\bSigma',j,\ell} \right\rVert_{\psi_2}\lesssim \sqrt{{\rm Var}(W_{\bSigma',j,\ell})}\le \lVert \bSigma_0^{1/2}(\bSigma'^{-1}-\bSigma_0^{-1})\rVert  \lVert \bX_j\rVert\lesssim \lVert \bSigma' -\bSigma_0  \rVert_F \lVert \bX\rVert_\circ.
	\end{align*} 
	We also obtain that
	\begin{align*}
	&\int_0^{c_5\lVert \bX\rVert_\circ{\rm diam}_j({\cal H}_n)}\sqrt{\log N\left(\frac{\epsilon}{2c_5\lVert \bX\rVert_\circ},{\cal H}_n ,\lVert\cdot\rVert_F\right)}d\epsilon\\
	&\quad\le \int_0^{2c_5\sqrt{M_1}\lVert \bX\rVert_\circ \epsilon_n }\sqrt{d^2 \log\left(\frac{6c_5\sqrt{M_1}\lVert\bX\rVert_\circ\epsilon_n}{\epsilon}\right)  }d\epsilon\\
	&\quad=  12 c_5\sqrt{M_1} \lVert \bX\rVert_\circ d \epsilon_n  \int_{\sqrt{\log 3}}^\infty t^2 e^{-t^2} dt.
	\end{align*} 
	Since the integral term on the rightmost side of the last display is bounded, we finally verify that $ \left\lVert  \sup_{\bSigma\in{\cal H}_n} | W_{\bSigma,j,\ell} | \right\rVert_{\psi_2}\lesssim \lVert \bX\rVert_\circ d\epsilon_n $ for every $j$ and $\ell$. Putting everything together, \eqref{eqn:upper22} is bounded by a multiple of $\epsilon_n^2\sqrt{p_{\max} n d^3 s^\star \log G} $ which goes to zero by the assumption. We finally verify \eqref{eqn:diffll}, and hence the proof is complete.
\end{proof}

\begin{proof}[Proof of Theorem \ref{thm-selection}]
	We only need to show that
	\begin{align*}
	\sup_{\substack{
			\bbeta_0 \in  \{\mathcal{B}_0 : s_0\le s_n \delta_n(s_0)\le \eta_n,\\ \phi_{\ell_{2,1}}(s_0+M_2s^\star) > c \} ,\,
			\bSigma_0 \in \mathcal{H}_0
	}}
	\mathbb{E}_0\Pi(\bbeta:S_{\bbeta,1}\supset S_{0,1} , \dots, S_{\bbeta,d}\supset S_{0,d} , S_{\bbeta} \ne S_0
	|Y_1, \dots, Y_n) 
	\rightarrow 0.
	\end{align*}
	Then the theorem follows by the Beta-min condition.
	Our proof is similar to the proof of Theorem 4 of \citet{cast15}.
	
	Let 
	$\mathcal{S}_n = 
	\{S: s \leq M_2 s^\star, S_1\supset S_{0,1} , \dots, S_d\supset S_{0,d} , S \ne S_0 \}$
	and 
	$\bGamma_{S} = 
	\sum_{i=1}^n \tilde\bX_{i,S} \bSigma_0^{-1} \tilde\bX_{i,S}'$.
	By Theorem~\ref{thm-bvm},
	it suffices to 
	show that
	$\Pi^\infty(\bbeta: S_{\bbeta} \in {\cal S}_n|Y_1, \dots, Y_n)
	\rightarrow 0$ in probability.
	By 
	(\ref{eqn-post-mixture-infty}), 
	we obtain that $
	\Pi^\infty(\bbeta: S_{\bbeta} \in \mathcal{S}_n 
	| Y_1, \dots, Y_n) \leq 
	\sum_{S \in \mathcal{S}_n}
	{w_S^\infty}/{ w_{S_0}^\infty}$ 
	which can be bounded by 
	\begin{align}
	&\sum_{\bar s=s_0+1}^{M_2s^\star}\Bigg\{
	\frac{
		\pi(\bar s){{Gd \choose s_0}}
		{{Gd - s_0} \choose {\bar s - s_0}}
	}
	{\pi(s_0){Gd \choose \bar s}}
	\max_{S\in {\cal S}_n : s=\bar s}\Bigg[
	\Bigg(
	\prod_{k=1}^d \prod_{j \in S_k} 
	\bigg(\frac{\lambda_j\sqrt{2\pi}}{a_j}\bigg)^{p_j}
	\Bigg)
	\left(\frac{\det \bGamma_{S_0}}{\det \bGamma_{S}}\right)^{1/2}
	\nonumber \\
	& \qquad \qquad \qquad \times 
	\exp\left(
	\frac{1}{2} \sum_{i=1}^n \|\bSigma_0^{-1/2}\tilde\bX_{i,S}'
	\hat\beta^\star_S\|^2 - 
	\frac{1}{2} \sum_{i=1}^n \|\bSigma_0^{-1/2}\tilde\bX_{i,S_0}'
	\hat\beta^\star_{S_0}\|^2 
	\right)\Bigg]\Bigg\}.
	\label{eqn-pthm4.4-1}
	\end{align}
	To bound further, we bound each factor in the above expression. 
	
	The interlacing theorem applied to $\bGamma_S$ and its principal submatrix $\bGamma_{S_0}$ gives $\eig_m(\bGamma_{S_0}) \leq \eig_m(\bGamma_S)$, $m = 1, \dots, \sum_{k=1}^d \sum_{j \in S_{0,k}} p_j$, we have
	$$
	\det (\bGamma_{S_0})
	\leq \prod_{m}
	\eig_m (\bGamma_S) 
	\leq \left(\eig_1(\bGamma_S) \right)^{p_{S_0}-p_S} \det (\bGamma_S),
	$$
	so by (\ref{assump-res}),
	${\det (\bGamma_{S_0})}/{\det (\bGamma_S)}$ is bounded by 
	$(b_2^{-1} \phi_{\ell_2}(s) \|\bX\|_\circ )^{2(p_{S_0}-p_S)}$.
	
	The exponential term  $Q_S := \sum_{i=1}^n \|\bSigma_0^{-1/2} \tilde\bX_{i,S}' \hat\beta_S^\star\|^2 -  \sum_{i=1}^n \|\bSigma_0^{-1/2} \tilde\bX_{i,S_0}' \hat\beta_{S_0}^\star\|^2$ in (\ref{eqn-pthm4.4-1}) has a $\chi^2$-distribution with degree of freedom $p_{S_0}-p_S$.
	By Markov's inequality on the exponential moment, we have that for every $0< u < 1/2$ and $r>0$,
	\begin{align*}
	& \mathbb{P}_0
	\left(
	\max_{S\in {\cal S}_n:s=\bar s} Q_S \geq r(\bar s-s_0) (\log G \vee p_{\max}\log n)
	\right) \\
	&\quad\le \exp\Big( -ur(\bar s -s_0) (\log G \vee p_{\max}\log n)\Big) \mathbb{E}_0 \left(\max_{S\in {\cal S}_n:s=\bar s} e^{u Q_S}\right) \\
	&\quad\le N_{\bar s}\exp\Big( -ur(\bar s -s_0) (\log G \vee p_{\max}\log n)\Big) (1-2u)^{-(p_{S_0}-p_S)/2},
	\end{align*}
	where $N_{\bar s}=\binom{Gd-s_0}{\bar s -s_0}$ is the cardinality of the set $\{S\in {\cal S}_n:s=\bar s\}$.
	Since $N_{\bar s}\le (Gd)^{\bar s -s_0}$ and $d^2\log n \ll n$, we have that for some $c>0$,
	\begin{align*}
	& \mathbb{P}
	\Big( Q_S \geq r(\bar s-s_0) (\log G \vee p_{\max}\log n), \text{ for any $S\in {\cal S}_n$}
	\Big) \\
	& \quad \leq 
	\sum_{\bar s > s_0} 
	\exp\left( -ur(\bar s -s_0) (\log G \vee p_{\max}\log n) + \frac{3}{2}(\bar s -s_0) \log G+c (\bar s-s_0)p_{\max} \right)
	\end{align*}
	which goes to $0$ whenever $ur>3/2$. If $r>3$, this is ensured by choosing $u$ arbitrarily close to $1/2$. Thus with probability tending to $1$, 
	(\ref{eqn-pthm4.4-1}) is bounded by 
	\begin{align}
	\sum_{s = s_0+1}^{M_2s^\star} 
	\frac{A_1^{s-s_0} s^{s-s_0} }{ (G \vee n^{p_{\max}})^{A_4(s-s_0)}}
	\left(
	\frac{\max_{1\le k\le d} \lambda_k\sqrt{2\pi }}{b_2^{-1}\|\bX\|_\circ \phi_{\ell_2}(s)}
	\right)^{p_{\max}(\bar s - s_0)} \frac{1}{(G \vee n^{p_{\max}})^{ r(\bar s-s_0)/2}}.
	\label{eqn:bbb}
	\end{align}
	Under the small $\lambda$ regime, for every $S$ such that $s\le M_2s^\star\lesssim G^a$,
	\begin{align*}
	\left(\frac{\max_{1\le k\le d} \lambda_k\sqrt{2\pi }}{b_2^{-1}\|\bX\|_\circ \phi_{\ell_2}(s)}
	\right)^{p_{\max}(\bar s - s_0)}\le  
	\left(\frac{\max_{1\le k\le d} \lambda_k  \sqrt{2\pi M_2 s^\star  }}{b_2^{-1}\|\bX\|_\circ \phi_{\ell_{2,1}}(M_2s^\star)}
	\right)^{p_{\max}(\bar s - s_0)}\lesssim 1 .
	\end{align*}
	and hence \eqref{eqn:bbb} goes to $0$ if $a-A_4+r/2<0$. If $A_4>a+3/2$, this is ensured by choosing $r$ arbitrarily close to $3$.
	
\end{proof}

%%%%%%%%%%%%%%%%%%%%%%%%%%%%%%%%%%%%%%
%%%%%%%%%%%%%%%%%%%%%%%%%%%%%%%%%%%%%%
%%%%%%%%%%%%%%%%%%%%%%%%%%%%%%%%%%%%%%
\section{Supplement to ``Bayesian Linear Regression for Multivariate Responses Under Group Sparsity''}

The following lemma obtains the normalizing constant in the density proportional to $e^{-\lambda \|x\|}$, $x=(x_1,\ldots,x_m)\in \mathbb{R}^m$. 

\begin{lemma}
	\label{lemma-apx-const-a}
	For $
	a_m = \sqrt{\pi} 
	\left(
	{\Gamma(m+1)}/{\Gamma({m}/{2}+1)}
	\right)^{{1}/{m}}$,
	\begin{equation}
	\label{S-eqn-lemmaA.1-1}
	\int_{\mathbb{R}^m} \left(\frac{\lambda}{a_m}\right)^m
	\exp(-\lambda \|(x_1,\ldots,x_m)\|) dx_1\cdots dx_m = 1.
	\end{equation}
	Also as $m\to\infty$, $a_m \asymp m^{1/2}$.
	
	If 	$x$ is expressed in terms of the spherical polar coordinates by a
	radius $r$, 
	a base angle $\theta_{m-1} \in (0,2\pi)$, 
	and $m-2$ angles $\theta_1, \dots, \theta_{m-2}$
	ranging over $(-\pi/2,\pi/2)$,
	then the density of $r$ is given by 
	\begin{align}
	\label{S-eqn-lemmaA.1-2}
	f(r|\lambda) = 
	\frac{\lambda^m}{\Gamma(m)}
	r^{m-1} \exp(-\lambda r),
	\end{align}
	the gamma density with the shape parameter $m$
	and rate parameter $\lambda$.
\end{lemma}

\begin{proof}
	Applying the polar transformation $(x_1,\ldots,x_m) \mapsto (r,\theta_1,\ldots,\theta_{m-1})$, evaluating the Jacobian,
	and applying the results shown in Chapter 1.5.1 of \citet{scot15},
	the integral in (\ref{S-eqn-lemmaA.1-1}) equals to
	\begin{align}
	%\label{S-eqn-lemmaA.1-2}
	& \int_0^{2\pi} \int_{-\pi/2}^{\pi/2} \cdots
	\int_{-\pi/2}^{\pi/2}
	\int
	\left(\frac{\lambda}{a}\right)^m
	r^{m-1} e^{-\lambda r} 
	\prod_{i = 1}^{m-2} (\cos\theta_{m-i-1})^i 
	d r\,
	d \theta_1 \cdots d \theta_{m-2}\, d\theta_{m-1}
	\nonumber \\
	& \quad = 
	\int_0^\infty \frac{2\pi^{m/2}\lambda^m}
	{\Gamma({m}/{2}) a^m}
	r^{m-1} e^{-\lambda r}
	d r.
	\label{S-eqn-lemmaA.1-3}    
	\end{align}
	The second line of the last display is obtained by using the results in 
	Chapter 1.5.2 of \citet{scot15}.
	Since $\int_0^\infty r^{m-1} e^{-\lambda r} dr=\Gamma(m)/\lambda^m$, the value
	\begin{align*}
	a_m = \sqrt{\pi} \left(\frac{2\Gamma(m)}{\Gamma({m}/{2})}
	\right)^{{1}/{m}}
	= \sqrt{\pi}\left(\frac{\Gamma(m+1)}{\Gamma({m}/{2}+1)}
	\right)^{{1}/{m}}
	\end{align*}
	makes $(\lambda/a_m)^m \exp\left(-\lambda \|(x_1,\ldots,x_m)\|\right)$ a probability density function. 
	
	Now by Stirling's approximation to the gamma functions,
	we obtain that
	$$\displaystyle
	\frac{\sqrt{2}\pi}{e}\left(\frac{2m}{e}\right)^{1/2}
	\leq a_m \leq 
	\frac{e}{\sqrt{2}}\left(\frac{2m}{e}\right)^{1/2},$$
	which implies that $a_m \asymp m^{1/2}$. The relation 
	(\ref{S-eqn-lemmaA.1-2}) is evident from (\ref{S-eqn-lemmaA.1-3}). 
\end{proof}

\begin{theorem}
	\label{rate Wishart}
	Consider the setup of Theorem~\ref{thm-3.1} except that the prior on $\bSigma$ is given by the inverse-Wishart distribution $\bSigma^{-1} \sim \mathcal{W}_d(\kappa d^2,  \bPhi)$ such that $c^{-1}\le\eig_1(\bPhi)\le \eig_d(\bPhi)\le c$ for some constant $c>1$, where $\mathcal{W}_d(\nu,  \bPsi)$ stands for the Wishart distribution in dimension $d$ with $\nu$ degrees of freedom and scale matrix $\bPsi$. 
	Then for a sufficiently large $M'>0$,
	\begin{align}
		%\label{equ-dim}
	& \sup_{\bbeta_0 \in \mathcal{B}_0,
		\bSigma_0 \in \mathcal{H}_0}
	\mathbb{E}_0
	\Pi \Big(\bbeta: s_{\bbeta}
	\geq
M' \tilde s^\star
	\Big | Y_1, \dots, Y_n
	\Big) \rightarrow 0, 
	\label{S-eqn:wisdim}\\
	& \sup_{\bbeta_0 \in \mathcal{B}_0,
		\bSigma_0 \in \mathcal{H}_0}
	\mathbb{E}_0  \Pi\Big(
	\bbeta: 
	\|\bX(\bbeta - \bbeta_0)\|_F^2
	\geq M' n\tilde\epsilon_n^2 \Big| Y_1, \dots, Y_n
	\Big) \rightarrow 0, 
	\label{S-eqn-post-beta}\\
	& \sup_{\bbeta_0 \in \mathcal{B}_0,
		\bSigma_0 \in \mathcal{H}_0}
	\mathbb{E}_0  \Pi\Big(
	\bSigma:
	\|\bSigma - \bSigma_0\|_F^2
	\geq M' \tilde\epsilon_n^2 \Big| Y_1, \dots, Y_n
	\Big) \rightarrow 0, 
	\label{Wishart-post-sigma}
	\end{align}
	where
	\begin{align}
	\label{S-eqn:sstartilde}
	\tilde s^\star &= \max\left\{s_0,  \frac{d^3\log n}{\log G\vee p_{\max}\log n}\right\},\\
	\label{Wishart-post-rate}
	\tilde\epsilon_n &= \max \left\{
	\sqrt{\frac{s_{0}  \log G}{n}} 
	,
	\sqrt{\frac{s_0 p_{\max}  \log n}{n}}
	,
	\sqrt{\frac{d^3\log n}{n}} \right\}.
	\end{align}
\end{theorem}

\begin{remark}\rm 
	\label{selection wishart}
	Once $\bbeta$ and $\bSigma$ are confined in small neighborhoods around the true values, the distributional approximation in Theorem~\ref{thm-bvm} and the selection consistency in Theorem~\ref{thm-selection} remain valid with the revised rate $\tilde\epsilon_n$ given in \eqref{Wishart-post-rate}. This can be shown by imitating the proofs of these theorems with the inverse Wishart prior on $\bSigma$, as the proofs of these results do not require a specific prior. 
\end{remark}

To prove the theorem, we need the following lemma giving estimates on the distribution of eigenvalues of a Wishart matrix. 

\begin{lemma}
	\label{lemma-A4}
	If $\bSigma^{-1} \sim \mathcal{W}_d(\nu,  \bPsi)$, 
	where
	$\nu \ge d$ is an integer, 
	$0 < \rho_1 < \cdots < \rho_d$ are its eigenvalues, 
	then 
	for $t_1 > \nu d$, $t_2 > 0$,  $0 \leq t_3 \leq 1$,
	and $0 \leq a_1 \leq \dots \leq a_d$,
	\begin{align}
	\mathbb{P}\left( \rho_d \geq t_1 \|\bPsi\|\right) 
	&\le 
	\left(
	\frac{t_1}{\nu d}
	\right)^{\nu d/2}
	\exp(\nu d/2- t_1/2),
	\label{S-eqn-lemmaA4-1}\\
\mathbb{P}\left( \rho_1 \leq t_2\right) 
&\le  \left(\frac{\nu+d}{2e}\right)^{d(\nu+d)/2} \frac{(e(\nu+d)/\sqrt{\pi})^d}{2^{(\nu+d+1)/2}}
t_2^{(\nu -d - 1)/2} \nonumber\\
&~~\times\left(\det(\bPsi)\right)^{-\nu /2}\|\bPsi\|^{(d-1)(\nu+1)/2},
\label{S-eqn-lemmaA4-2} \\
	\mathbb{P}\left(\bigcap_{k=1}^d 
	\{a_k \leq \rho_k \leq a_k(1+t_3)\}\right) 
& \ge  
	\left(
	\frac{a_1 t_3 e^2\nu}{8\sqrt{\pi}}
	\right)^{-d}
	\left(
	\frac{2\nu d}{ea_1t_3}
	\right)^{-\nu d/2}
	\left(
	\frac{d}{2e}
	\right)^{-d^2/2}\nonumber\\
	& 	~~\times\left(\det(\bPsi)\right)^{-\nu /2}
	\exp\left(- \frac{a_1(1+t_3)\Tr(\bPsi^{-1})}{2}\right).
	\label{S-eqn-lemmaA4-3}
	\end{align}
\end{lemma}

\begin{remark} \rm 
	\label{special Wishart}
	To control the sequences appearing in the estimates of prior probabilities in Lemma~\ref{lemma-A4} such that explicit growth estimates can be obtained for use in the rate theorem, the degrees of freedom of the Wishart prior on $\bSigma^{-1}$ needs to be taken approximately proportional to the dimension $d^2$. By choosing $\nu$ to be the integer part of $\kappa d^2$ for some constant $\kappa \ge 1$, 
	the estimates in Lemma~\ref{lemma-A4} simplify to 
	\begin{align*}
	\mathbb{P}\left( \rho_d \geq t_1 \|\bPsi\|\right) 
	&\le 
	\left(
	b_1{t_1}/{d^3}
	\right)^{b_2d^3}
	\exp(b_3 d^3- t_1/2),
	\nonumber%\label{S-eqn-lemmaA4-1.1}
	\\
	\mathbb{P}\left( \rho_1 \leq t_2\right) 
	&\le 
	\left({b_4d^2}\right)^{b_5d^3}
	t_2^{b_6d^2}
\left(\det(\bPsi)\right)^{-\kappa d^2 /2}\|\bPsi\|^{b_7d^3}
	\nonumber%\label{S-eqn-lemmaA4-2.1} 
	\\
	 \mathbb{P}\left(\bigcap_{k=1}^d 
	\{a_k \leq \rho_k \leq a_k(1+t_3)\}\right) 
	& \ge 
	(b_8a_1 t_3 d^2)^{-d}(b_9d^3/(a_1t_3))^{-b_{10}d^3}
	(b_{11}d)^{-d^2/2} 
	\nonumber
	\\
&  \quad\times 	\left(\det(\bPsi)\right)^{-\kappa d^2/2}
	\exp\left(-{a_1(1+t_3)\Tr(\bPsi^{-1})}/{2}\right),
	%\label{S-eqn-lemmaA4-3.1}
	\end{align*}
	for some constants $b_1,\dots, b_{11}>0$.
\end{remark}

\begin{remark}\rm
	\label{selection wishart2}
	In Theorem~\ref{rate Wishart}, the degree of freedom $\nu$ is chosen to grow in proportion to $d^2$.
	This choice makes the first two terms of $\tilde\epsilon_n$ the same as those in $\epsilon_n$, but induces slightly increased $\tilde s^\star$ in \eqref{S-eqn:sstartilde} compared to $s^\star$ as the prior concentration decreases.
	Instead, one may choose $\nu$ that is proportional to $d$. Then the prior concentration stays and the assertion \eqref{S-eqn:wisdim} holds with $s^\star$ instead of $\tilde s^\star$, so the same dimension recovery result is obtained as Lemma~\ref{lemma-dim}. However, this significantly weakens the rate to $ \max \{\sqrt{ (ds_{0}  \log G)/n} ,\sqrt{ (ds_0 p_{\max}  \log n)/n},\sqrt{ (d^3\log n )/n} \}$, because a weaker right-tail decay of the largest eigenvalue of $\bSigma$ necessitates the use of a larger sieve, which increases the entropy.
\end{remark}

\begin{proof}[Proof of Lemma~\ref{lemma-A4}]
	If the dimension $d$ remains bounded, the result is already given in Lemma~9.16 of \citet{ghos17}. For $d\to \infty$, the dependence of the constants on $d$ must be explicitly identified. Below we carefully estimate the normalizing constant, assuming that $d$ is sufficiently large. 
	
	To prove (\ref{S-eqn-lemmaA4-1}), consider the random matrix  
	 $\bOmega = \bPsi^{1/2} \bSigma \bPsi^{1/2}$. Observe that 
	then $\bOmega^{-1} \sim \mathcal{W}_d(\nu, \bI_d)$. 
	Since $\rho_d = \|\bSigma^{-1}\| = \|\bPsi^{1/2}\bOmega^{-1}\bPsi^{1/2}\|
	\leq \|\bPsi\|\|\bOmega^{-1}\|$,
	we have that 
	$$\mathbb{P}(\rho_d \geq t_1\|\bPsi\|) 
	\leq \mathbb{P}(\|\bOmega^{-1}\| \geq t_1) 
	\leq \mathbb{P}(\Tr(\bOmega^{-1}) \geq t_1)$$ 
	and $\Tr(\bOmega^{-1}) \sim \chi_{\nu d}^2$.
	Then apply the Chernoff bound for a $\chi^2$-distribution,
	we obtain (\ref{S-eqn-lemmaA4-1}).
	
	To prove (\ref{S-eqn-lemmaA4-2}) and (\ref{S-eqn-lemmaA4-3}),
	we need estimates for the multivariate gamma function from both sides.	
	By Stirling's approximation to the gamma functions,
	$e(n/e)^n \leq \Gamma(n+1) \leq en(n/e)^n$. Thus  
	we have  
	\begin{align*}
	\Gamma_d(\nu/2)
	& = \pi^{d(d-1)/4}
	\prod_{k=1}^d \Gamma
	\left(
	\frac{\nu+1-k}{2}
	\right)\\
&	\leq 
	\pi^{d(d-1)/4}
	\left( \Gamma
	\left(
	{\nu}/{2} + 1
	\right)
	\right)^d
	\\
	& 
	\leq 
	\pi^{d(d-1)/4}
	e^{-\nu d/2 + d}
	\left({\nu}/{2}\right)^{(\nu/2 + 1)d}
	\end{align*}
	and since $\nu \geq d$,
	\begin{align*}
	\Gamma_d(\nu/2)
	& 
	\geq 
	\pi^{d(d-1)/4}
	(\Gamma(1/2))^d
	= 
	\pi^{d(d-1)/4 + d/2}.
	\end{align*}

	To prove (\ref{S-eqn-lemmaA4-3}),
	we need the following three inequalities: 
	\begin{enumerate} 
		\item $\displaystyle 
	%\begin{align*}
	\prod_{1\le k < k'\le d} (\rho_{k'} - \rho_k) 
	\leq \prod_{1\le k < k'\le d} \rho_{k'} = \prod_{k=2}^d \rho_k^{k-1}; 
	%\end{align*}
$	
	\item $\displaystyle 
%	\begin{align*}
	\exp\left(-\frac{\Tr(\bPsi^{-1}\bP \bD^{-1} \boldsymbol P')}{2}\right) 
	 \leq
	\exp\left(-\frac{ \Tr(\bP \bD^{-1} \bP')}{2\|\bPsi\|}\right) 
	 = \exp\left(-\sum_{k=1}^d \frac{\rho_k}{2\|\bPsi\|} \right);
%	\end{align*}
$
\item $\displaystyle
%\begin{align*}
\frac{\pi^{d^2/2} 2^{-d\nu/2}
	\left(\det(\bPsi)\right)^{-\nu/2}}
{\Gamma_d(d/2) \Gamma_d(\nu/2)}
\leq \pi^{-d/2} 2^{-d\nu/2}\left(\det(\bPsi)\right)^{-\nu/2},
%\end{align*}
$
which is a consequence of the lower bound for the multivariate gamma function. 
\end{enumerate}
	
	Note that $\bD^{-1}={\rm diag}(\rho_1,\dots,\rho_d)$.
	Then by plugging-in the above three upper bounds in the expression for the joint density of the eigenvalues of a Wishart matrix (see, e.g., equation (9.6) of \citet{ghos17}) and integrating, the marginal density of $\rho_1$ is bounded by 
	\begin{align*}
	& \pi^{-d/2} 2^{-d\nu/2}\left(\det(\bPsi)\right)^{-\nu/2}
	\rho_1^{(\nu - 1 - d)/2} e^{-\rho_1/(2\|\bPsi\|)}\\
	&  \quad  \times \prod_{k=2}^d 
	\int_{0}^{\infty}
	\rho_k^{(\nu-1-d)/2+k-1} \exp(-\rho_k/(2\|\bPsi\|)) d \rho_k.
	\end{align*}
	Each integral function in the last display equals to 
	$ \Gamma \left(
	\frac{\nu-d-1}{2} + k
	\right)
	\left(
	2\|\bPsi\|
	\right)^{(\nu - d -1)/2 + k}$. 
Now applying the upper bound of the gamma function, we obtain 
	\begin{align*}
	\Gamma \left(
	\frac{\nu-d-1}{2} + k
	\right)
	\le
	\Gamma \left(
	\frac{\nu+d}{2} + 1
	\right)
	\leq 
	\left(\frac{\nu + d}{2}\right)^{(\nu+d)/2 + 1} 
	e^{1 -{\nu + d}/{2}}.
	\end{align*}
	With $\rho_1 \leq t_2$, the marginal density of $\rho_1$ can be further bounded above by
	\begin{align*}
	\pi^{-d/2}2^{-\nu d /2} 
	(\det(\bPsi))^{-\nu /2} 
	\left(\frac{\nu + d}{2}\right)^{d(\nu+d)/2 + d}
	t_2^{(\nu-d-1)/2} e^{-d(\nu + d)/2 + d} 
	\left(
	2\|\bPsi\|
	\right)^{(d-1)(\nu+1)/2},
	\end{align*}
	as $\sum_{k=2}^d ((\nu-d-1)/2+k)=(d-1)(\nu+1)/2$ if $d\ge 2$.
	The last display equals to the upper bound of (\ref{S-eqn-lemmaA4-2}).

	To prove (\ref{S-eqn-lemmaA4-3}),
	let 
	$I_k = \{a_k(1 + (k-1/2)t_3/d), a_k(1 + kt_3/d)\}$ for each $k \in \{1, \dots, d\}$. Then $\rho_k \in I_k$ implies that
	$\rho_k \in [a_k, a_k(1+t_3)]$.
	Therefore integrating the expression for the joint density of the eigenvalues and using the estimates of the normalizing constants given above, we have
	\begin{align}
	& \mathbb{P}\left(\bigcap_{k=1}^d 
	\{\bSigma: a_k \leq \rho_k \leq a_k(1+t_3)\}\right)
	\nonumber \\
	& \quad \geq
	\frac{\pi^{d^2/2} 2^{-d\nu/2}
		\left(\det(\bPsi)\right)^{-\nu/2}}
	{\Gamma_d(d/2) \Gamma_d(\nu/2)}
	\int_{I_d} \dots \int_{I_1}\bigg\{
	\prod_{k=1}^d \rho_k^{(\nu-d-1)/2} 
	\prod_{k<k'}^d (\rho_{k'} - \rho_{k})
	\nonumber	\\
	& \qquad\qquad \qquad\qquad \qquad\qquad \qquad \times
	\int_{\mathscr{O}(d)} 
	\exp\left( - \frac{1}{2}  \text{Tr}(
	\bPsi^{-1}
	\bP \bD^{-1}  \bP' )\right)\bigg\}
	d \bP\, d \rho_1 \cdots d \rho_d
	\nonumber \\
	& \quad  \geq
	\frac{\pi^{d^2/2} 2^{-d\nu/2}
		\left(\det(\bPsi)\right)^{-\nu/2}}
	{\Gamma_d(d/2) \Gamma_d(\nu/2)}
	\left(\frac{a_1 t_3}{2d}\right)^{(\nu-2)d/2} 
	\exp\left(- \frac{a_d(1+t_3)}{2}
	\text{Tr}( \bPsi^{-1})
	\right).
	\label{S-eqn-lemmaA4-4}
	\end{align}
	The lower bound in the third line of the last display is obtained by noticing that for $k' > k$, $\rho_{k'} - \rho_{k} \geq {a_1 t_3}/{(2d)}$,
	$-\bD^{-1} \ge -\rho_d \bI_d > - a_d(1+t_3) \bI_d$,
	and $ \bP  \bP' =  \bI_d$. 
	Now we plug the upper bound for the multivariate gamma function in (\ref{S-eqn-lemmaA4-4}) to obtain the lower bound in 
	(\ref{S-eqn-lemmaA4-3}).
\end{proof}

\begin{proof}[Proof of Theorem~\ref{rate Wishart}]
	It only suffices to obtain estimates of prior concentration and define an appropriate sieve for this prior such that the complement has exponentially small prior probability. The proof is very similar to that of Theorem~\ref{thm-3.1} employing the same overall strategy, except when estimates regarding the prior concentration of the covariance matrix are involved. The estimates of the prior mass outside the sieve and that of the entropy of the sieve must be obtained afresh since a different sieve is used.

Since the negative logarithm of the average Kullback-Leibler neighborhood of size $\tilde\epsilon_n^2$ should be controlled, the probabilities of the sets $\{\bbeta: \sum_{k=1}^d\lVert\beta_k-\beta_{0,k}\rVert_{2,1}
\leq c \tilde r_n\}$ and $\{\bSigma: \|\bSigma^{*-1}-\bI\|\le \tilde\epsilon_n\}$ need to be obtained, where $\tilde r_n = 
\sqrt{{n\tilde\epsilon_n^2}}/{ \|\bX\|_\circ}$ and $\bSigma^*= \bSigma_0^{-1/2}\bSigma\bSigma_0^{-1/2}$ as before.  
For the former, it is easy to see that the prior concentration of $\bbeta$ is bounded by a constant multiple of $n\tilde\epsilon_n^2$, in view of the proof of Lemma~\ref{lemma-3.1}.
The condition for the latter clearly holds if all eigenvalues of $\bSigma^{*-1}$ lie between $1$ and $1+d^{-1/2}\tilde\epsilon_n$. In view of the third assertion in Remark~\ref{special Wishart}, the prior probability of this event is at least 	
	$$
	-d\log (b_8d^{3/2}\tilde\epsilon_n)-b_{10}d^3 \log\left(\frac{b_9d^{7/2}}{\tilde\epsilon_n}\right)-\frac{d^2}{2}\log(b_{11} d)-\frac{\kappa d^3}{2} \log\lVert\bPsi\rVert-\frac{d+d^{1/2}\tilde\epsilon_n}{2}\log\lVert\bPsi^{-1}\rVert,
$$
which is bounded below by a constant multiple of $-d^3\log n$.
Thus, the estimate for the prior concentration is controlled. Then using the same techniques in the proof of Lemma~\ref{lemma-dim}, dimension recovery is still valid with $s^\star$ replaced by $\tilde s^\star$, which verifies \eqref{S-eqn:wisdim}.

	Next, to prove \eqref{S-eqn-post-beta} and \eqref{Wishart-post-sigma} define the sieve
	\begin{align*}
	\tilde{\mathcal{F}}_n = 
	\Big\{
	(\bbeta,\bSigma)\in{\cal B}_n\times{\cal H}: &\   s\le M' \tilde s^\star,\
	\max_{\substack{1\leq j\leq G\\ 1 \leq k \leq d}}
	\|\beta_{jk}\|
	\le  H_n, \\
	&
	\exp(-M n\tilde\epsilon_n^2/d^2) < \eig_1({\bSigma}^{-1}), \
	\eig_d({\bSigma}^{-1}) 
	\le n
	\Big\},
	\end{align*}
	for a sufficiently large $M>0$,
	where $	H_n = {p_{\max}n}/\underline\lambda$.
	Recall that the expression of $\underline{\lambda}$ is shown in (\ref{eqn-assump-2}).

	We shall verify that 
	\begin{align}
	\Pi( ({\cal B}_n\times {\cal H})\setminus\tilde{\mathcal{F}}_n) 
	\leq \exp\left(-(1+C_1)n \tilde\epsilon_n^2\right),
	\label{S-eqn-pthm3.1-2}
	\end{align}
	as long as $M$ is chosen sufficiently large.

Following \eqref{eqn-pthm3.1-11}, it suffices to bound only 
	the terms $\Pi\left(\eig_1({\bSigma}^{-1}) \leq \exp(-M n\tilde\epsilon_n^2/d^2)\right)$ and  
	$\Pi\left(\eig_d({\bSigma}^{-1}) \geq n\right)$ as the rest is unchanged.
	
	In view of Remark~\ref{special Wishart}, these terms are  bounded above by
	\begin{align*}
	&\exp \left(c_1 d^3\log n - c_2 M n\tilde\epsilon_n^2\right) + 
	\exp\left(c_3d^3\log n - c_4 n\right),
	\end{align*}
	where $c_{1}, \dots, c_{4}$ are positive constants.
Thus if $M$ is chosen sufficiently large, we have \eqref{S-eqn-pthm3.1-2}.
	
	To complete the proof of \eqref{eqn-pthm3.1-3}, we need to show that $\log N_\ast \lesssim n\epsilon_n^2$, where $N_\ast$ is
	the number of pieces satisfying \eqref{eqn:smallpiece} needed to cover the sieve $\tilde{\mathcal F}_n$. It is easy to see that $\log N_\ast$ is bounded by
	\begin{align*}
	&\log N\Big(\frac{1} { \tilde s^\star \sqrt{p_{\max}n} \lVert \bX \rVert_\circ}, \Big\{\bbeta : s_{\bbeta}\le M' \tilde s^\star,	\max_{\substack{1\leq j\leq G \\ 1 \leq k \leq d}}
	\|\beta_{jk}\|
	\le  H_n\Big\}, \lVert \cdot \rVert_\infty \Big)\\
	&\quad+ \log N \Big( \frac{1}{n^2 d} ,\left\{ \bSigma: 	\exp(-M n\tilde\epsilon_n^2/d^2) < \eig_1({\bSigma}^{-1}), \
	\eig_d({\bSigma}^{-1}) 
	< n\right\}, \lVert\cdot\rVert\Big).
	\end{align*}
Similar to \eqref{eqn:entrobeta}, it can be easily verified that the estimate of the first term is bounded by a constant multiple of $n\tilde\epsilon_n^2$, so we only need to estimate 
	the second term, which is bounded by
	\begin{align*}	
	&\log N \left( \frac{1}{n^2 d} ,\left\{ \bSigma: 	\exp(-n\tilde\epsilon_n^2/d^2) \le \eig_1({\bSigma}^{-1}) \right\}, \lVert\cdot\rVert\right)\\
	&\quad \le 	\log N \left( \frac{1}{n^2 d} ,\left\{ \bSigma: 	\lVert\bSigma\rVert_F < \sqrt{d}\exp(M n\tilde\epsilon_n^2/d^2) \right\}, \lVert\cdot\rVert_F\right)\\
	&\quad \le {d^2} \log  \Big( n^2 d^{3/2} \exp(M n\tilde \epsilon_n^2/d^2) \Big). 
	\end{align*}
The last expression is easily seen to be bounded by a constant multiple of $n\tilde\epsilon_n^2$. 
	
\end{proof}

\bibliographystyle{chicago}
\bibliography{citation}

\end{document}